\documentclass[12pt]{report}
\usepackage[centertags]{amsmath}
\usepackage{amscd}
\usepackage{amsfonts}
\usepackage{amssymb}
\usepackage{amsthm}
\usepackage{newlfont}
\usepackage{xthesis} 
\usepackage{xtocinc} 
\usepackage{makeidx}
\usepackage{graphicx}

\usepackage[active]{srcltx}  
\RequirePackage{fancyhea}

\renewcommand{\baselinestretch}{1.6}

\newcommand{\Real}{\mathbb R}

\newcommand{\lc}{\mathcal{L}}
\newcommand{\g}{\mathcal{G}}

\newcommand{\Co}{\operatorname{Cone}}

\newcommand{\n}{{\mathcal{V}}}
\newcommand{\m}{{\mathcal{U}}}

\theoremstyle{plain}
\newtheorem{thm}{Theorem}[section]
\newtheorem{cor}[thm]{Corollary}
\newtheorem{lem}[thm]{Lemma}
\newtheorem{prop}[thm]{Proposition}
\theoremstyle{definition}
\newtheorem{defn}{Definition}[section]

\theoremstyle{remark}
\newtheorem{rem}{Remark}[section]
\theoremstyle{remark}
\newtheorem{exmp}{Example}[section]
\numberwithin{equation}{section}

\usepackage[cmtip, matrix, arrow]{xy}
\begin{document}



\dedicate{To Zahra, Ali and Reza}

\nolistoftables \nolistoffigures \phd

\copyrightyear{2004} \submitdate{April 2004} \convocation{}{2004}


\title{Differential Geometry of Relative Gerbes}

\author{Zohreh Shahbazi}

\supervisor{Eckhard Meinrenken}

\examiner{}

\firstreader{} \secondreader{}

{
\typeout{:?0000} 
\beforepreface
\typeout{:?1111} 
}
{
\typeout{Abstract}
\include{ABS}
}
\setcounter{page}{1} \pagenumbering{arabic} {
\typeout{Introduction}

\chapter{Introduction}

The motivation of this thesis is to develop pre-quantization of
quasi-Hamiltonian spaces with group-valued moment maps by
introducing the notion of \emph{relative gerbes} and addressing
their differential geometry.
\\ \indent Giraud~\cite{Gr} first introduced the concept
of gerbes in the early 1970s to study non-Abelian second
cohomology. Later, Brylinski~\cite{MR94b:57030} defined gerbes as
sheafs of groupoids with certain axioms, and discussed their
differential geometry. He proved that the group of equivalence
classes of gerbes gives a geometric realization of integral three
cohomology classes on manifolds. Through a more elementary
approach, Chatterjee and Hitchin \cite{D,MR2003f:53086} introduced
gerbes in terms of transition line bundles for a given cover of
the manifold. From this point of view, a gerbe is a one-degree-up
generalization of a line bundle, where the line bundle is
presented by transition maps. A notable example of a gerbe arises
as the obstruction for the existence of a lift of a principal
$G$-bundle to a central extension of the Lie group. Another
example is the associated gerbe
 of an oriented codimension three submanifold of an oriented
 manifold. The third example is what is called ``basic gerbe,'' which corresponds to the
 generator of the degree three integral cohomology of a compact,
 simple and simply connected Lie group. The basic gerbe over G is closely related to
 the basic central extension of the loop group, and it was constructed, from this point of
 view, by Brylinski~\cite{MR94b:57030}. Later,
 Gawedski-Reis~\cite{MR2003m:81222},
 for G=SU(n), and Meinrenken~\cite{EM}, in the general case, gave a finite-dimensional
 construction along with an explicit description of the gerbe
 connection.\\
\indent Quasi-Hamiltonian $G$-spaces $(M,\omega,\Psi)$ with
group-valued moment maps $\Psi:M\rightarrow G$ are introduced in
\cite{MR99k:58062}. The 2-form $\omega$ is not necessarily
non-degenerate. However, its kernel is controlled by the
\emph{minimal degeneracy} condition. Conjugacy classes of $G$,
with moment map the inclusion into $G$, are the main examples of
the quasi-Hamiltonian $G$-spaces. Another example is the space
$G^{2h}$, with moment map the product of Lie group commutators.
One can define the notion of reduction for group-valued moment
maps, and the reduced spaces $M//G=\Psi^{-1}(e)/G$ are symplectic.
The ``fusion product'' of two quasi-Hamiltonian $G$-spaces is a
quasi-Hamiltonian $G$-space, with moment map the pointwise product
of the two moment maps. This product
gives a ring structure to the set of quasi-Hamiltonian $G$-spaces.\\
\indent This thesis introduces the notion of \emph{relative
 gerbes} for smooth maps of manifolds, and discusses their
 differential geometry. The equivalence classes of relative gerbes
 are classified by the relative integral cohomology in degree
 three. Furthermore, by using the concept of relative gerbes, the
 pre-quantization of Lie group-valued moment maps is developed, and
 its equivalence with the pre-quantization of infinite-dimensional
 Hamiltonian loop group spaces is established.\\
 \indent The organization of this thesis is as follows. In Chapter
2, I discuss the relative (co)homology of a smooth map between two
manifolds. When the map is inclusion, the singular relative
(co)homology of the map coincides with the singular relative
(co)homology of the pair. Also, for a continuous map of
topological spaces, the relative (co)homology of the map is
isomorphic to the (co)homology of the mapping cone. In Chapter 3,
following the Chatterjee-Hitchin perspective on gerbes, I define
the notion of \emph{relative gerbe} for a smooth map $\Phi\in
C^{\infty}(M,N)$ between two manifolds $M$ and $N$ as a gerbe over
the target space together with a quasi-line bundle for the
pull-back gerbe. I prove that the group of equivalence classes of
relative gerbes can be characterized by
the integral degree three relative cohomology of the same map.\\
\indent Another objective of this thesis is to develop the
differential geometry of relative gerbes. More specifically, in
Chapter 4, the concepts of relative connection, relative
connection curvature, relative Cheeger-Simon differential
character, and relative holonomy are introduced. As well, I prove
that a given closed relative 3-form arises as a curvature of some
relative gerbe with connection if and only if the relative 3-form
is integral. I also prove that a relative gerbe with connection
for a smooth map $\Phi: M\to N$ generates a \emph{relative} line
bundle with connection for the corresponding map of loop paces,
$L\Phi:\,LM\to LN$. In addition, I prove that the group of
equivalence classes of relative gerbes with connection gives a
geometric realization of degree two relative Deligne cohomology.\\
\indent In Chapter 5, I give an explicit construction of the basic
gerbe for $G=SU(n)$, and of suitable multiples of the basic gerbe
for the other Lie groups. As well, I show that the construction of
the basic gerbe over $SU(n)$ is equivalent to the construction of
the basic gerbe in Gawedzki-Reis~\cite{MR2003m:81222}. In this
chapter, I also construct a relative gerbe for the map
$\operatorname{Hol}:\,\mathcal{A}_G(S^1)\rightarrow G$, where
$\mathcal{A}_G(S^1)$ is the affine space of connections on the
trivial bundle $S^1\times G$. Inspired by the pre-quantization of
Hamiltonian $G$-manifolds, the other objective of this thesis is
to construct a method to pre-quantize the quasi-Hamiltonian
$G$-spaces with group-valued moment maps. To achieve this, I use
the premise of \emph{relative gerbes}, developed in this thesis.
Let $G$ be a compact, simple and simply connected Lie group and
$(M,\omega,\Psi)$ be a quasi-Hamiltonian $G$-space. One of the
axioms of a group-valued moment map $(M,\omega,\Psi)$ is
$d\omega=\Psi^*\eta$, where $\eta$ is the canonical 3-form. This
means that the pair $(\omega,\eta)$ is a cocycle for the relative
de Rham complex for the moment map $\Psi$. A pre-quantization
 of a quasi-Hamiltonian $G$-space with $G$-valued moment map is defined in Chapter 5 by a
 relative gerbe with connection with relative curvature $(\omega,\eta)\in\Omega^3(\Psi)$.
 Based on the results of Chapter 4, $(M,\omega,\Psi)$
 is pre-quantizable if and only if the relative 3-form
 $(\omega,\eta)$ is integral.
  I prove that, given two pre-quantizable quasi-Hamiltonian $G$-spaces, their fusion product
  is again pre-quantizable.
 \\ \indent In Chapter 5, the pre-quantization conditions for the examples of
 quasi-Hamiltonian $G$-spaces, described previously, are examined. I show that a given
  conjugacy class $\mathcal{C}=G\cdot\exp(\xi)$ of $G$ is  pre-quantizable when $\xi\in
  \Lambda^*$, where $\Lambda^*$ is the weight lattice. I also illustrate that
  $G^{2h}$ has a unique pre-quantizaion, which enables one to construct
 a finite dimensional pre-quantum line bundle for the moduli space of flat
 connections of a closed oriented surface of genus $h$.

 Further, recall that there is a one-to-one correspondence between quasi-Hamiltonian $G$-spaces and
 infinite dimensional loop group spaces~\cite{MR99k:58062}.
 Extending this correspondence, in Chapter 5, I prove that pre-quantization of a
 quasi-Hamiltonian $G$-spaces with group-valued moment map coincides with the
 pre-quantization of the corresponding Hamiltonian loop group space.

 As a continuation of this thesis, one can consider the
 ``quantization'' of quasi-Hamiltonian $G$-spaces. It is expected that this advancement should
 involve the \emph{relative} notion of twisted K-theory, to be developed.

}


\chapter
{Relative Homology/Cohomology}


\section{Algebraic Mapping Cone for Chain Complexes }
\begin{defn}Let $f_{\bullet}:X_\bullet\longrightarrow Y_\bullet$
be a chain map between chain complexes over R where R is a
commutative ring. The algebraic mapping cone of
$f$~\cite{MR2002f:55001} is defined as a chain complex
$\operatorname{Cone}_{\bullet}(f)$ where
\[\operatorname{Cone}_n(f)=X_{n-1}\oplus Y_{n}\]with the
differential\[\partial(\theta,\eta)=(\partial\theta,f(\theta)-\partial
\eta).\]Since $\partial^2=0$, we can consider the homology of this
chain complex. Define relative homology of $f_{\bullet}$ to be
$$H_n(f):=H_n(\operatorname{Cone}_{\bullet}(f)).$$
\end{defn}
The short exact sequence of chain complexes\[0\rightarrow
Y_{n}\overset{j}\rightarrow
\operatorname{Cone}_n(f)\overset{k}\rightarrow X_{n-1}\rightarrow
0\] where $j(\beta)=(0,\beta)$ and $k(\alpha, \beta)=\alpha$ gives
a long exact sequence in
homology\begin{equation}\label{eq0}\begin{split}\cdots \rightarrow
H_{n}(Y)\overset{j}\rightarrow H_{n}(f)\overset{k}\rightarrow
H_{n-1}(X)\overset{\delta}\rightarrow H_{n-1}(Y) \rightarrow
\cdots\end{split}\end{equation}where $\delta$ is the connecting
homomorphism.
\begin{lem}The connecting homomorphism $\delta$ is given by
$\delta[\gamma]=[f(\gamma)]$ for $\gamma\in
X_{n-1}$.\end{lem}\begin{proof} For $\gamma\in X_{n-1}$, we have
$k(\gamma,0)=\gamma$. The short exact sequence of chain complexes
gives an element $\gamma'\in Y_{n-1}$ such that
$j(\gamma')=\partial(\gamma,0)=(\partial\gamma,f(\gamma))$.
$\delta$ is defined by $\delta[\gamma]=[\gamma']$. But by
definition of $j$, $j(\gamma')=(0,\gamma')$. Therefore
$f(\gamma)=\gamma'$. This shows
$\delta[\gamma]=[f(\gamma)]$.\end{proof}\begin{defn}We call a
chain map $f_{\bullet}:X_{\bullet}\rightarrow Y_{\bullet}$ a
quasi-isomorphism if it induces isomorphism in cohomology, i.e.,
$H_{\bullet}(X)\overset{\cong}\rightarrow
H_{\bullet}(Y)$.\end{defn}\begin{cor}
$f_{\bullet}:X_{\bullet}\rightarrow Y_{\bullet}$ is a
quasi-isomorphism if and only if
$H_{\bullet}(f)=0$.\end{cor}\begin{proof} $f$ is a
quasi-isomorphism, if and only if the connecting homomorphism in
the long exact sequence \ref{eq0}is an isomorphism.
\end{proof} \begin{defn}A homotopy operator between two chain
complexes
 $f,g:X_{\bullet}\rightarrow Y_{\bullet}$ is a linear map
 $h:X_{\bullet}\rightarrow Y_{\bullet+1}$
 such that \[h\partial+\partial h=f-g.\quad(\star)\]In that case, $f$ and $g$ are
 called chain homotopic and we denote it by $f\simeq
 g$.\newline \indent Two chain maps $f:X_{\bullet}\rightarrow
 Y_{\bullet}$ and $g:Y_{\bullet}\rightarrow
 X_{\bullet}$ are called homotopy inverse if $g\circ
 f\simeq id_X$ and $f\circ
 g\simeq id_Y$ are both homotopic to the identity. If
  $f:X_{\bullet}\rightarrow Y_{\bullet}$ admits a homotopy inverse, it is
  called a homotopy equivalence. In particular every homotopy equivalence is a
  quasi-isomorphism. \end{defn}
 \begin{prop}Any homotopy between chain maps $f,g:X_\bullet\rightarrow Y_\bullet$
 induces an isomorphism of chain complexes $\operatorname{Cone}(f)_{\bullet}$ and
 $\operatorname{Cone}(g)_{\bullet}$.
\end{prop}
  \begin{proof}Given a homotopy operator $h$ satisfying $(\star\,)$, define a map $F:\operatorname{Cone}_\bullet(f)\rightarrow
  \operatorname{Cone}_\bullet(g)$ by
  \[F(\alpha,\beta)=(\alpha,-h(\alpha)+\beta).\]

 Since\[
 \partial F(\alpha,\beta)=(\partial\alpha,g(\alpha)+\partial
 h(\alpha)+\partial\beta)=(\partial\alpha,f(\alpha)-h\partial(\alpha)+\partial\beta)
 =F\partial(\alpha,\beta)\]
 $F$ is a chain map and its inverse map is $F^{-1}(\alpha,\beta)=(\alpha,h(\alpha)+\beta)$.

  \end{proof}
\begin{lem}\label{1}Let \\*\[\begin{CD}
  0 @>>> X_{\bullet}  @>>>  Y_{\bullet} @>>> Z_{\bullet} @>>> 0  \\
   @.  @VVV  @VVV  @VVV  @.\\
   0 @>>>\widetilde{X}_{\bullet}@>>>\widetilde{Y}_{\bullet}@>>>\widetilde{Z}_{\bullet}
   @>>>0
   \end{CD}\]\\*be a commutative diagram of chain maps with exact
   rows. If two of vertical maps are quasi-isomorphism, then so is
   the third.\end{lem}
   \begin{proof}The statement follows from the 5-Lemma applied to
   the corresponding diagram in homology, \\*\[\begin{CD}
  \cdots @>>> H_{\bullet}(X)  @>>>  H_{\bullet}(Y) @>>> H_{\bullet}(Z) @>>> H_{\bullet-1}(X) @>>> \cdots \\
   @.  @VVV  @VVV  @VVV @VVV @.\\
   \cdots @>>>H_{\bullet}(\widetilde{X})@>>>H_{\bullet}(\widetilde{Y})@>>>H_{\bullet}(\widetilde{Z})
   @>>>H_{\bullet-1}(\widetilde{X})@>>>\cdots
   \end{CD}\]\end{proof}
  \begin{prop}\label{2}Suppose we have the following commutative diagram of chain
  maps,\\*\[\begin{CD}
   X_{\bullet}  @>f_{\bullet}>>     Y_{\bullet}\\
   @V\Phi_{\bullet}VV   @V\Psi_{\bullet}VV\\
   \widetilde{X}_{\bullet}@>\widetilde{f}_{\bullet}>>\widetilde{Y}_{\bullet}
   \end{CD}\]\\*such that $\Phi$ and $\Psi$ are quasi-isomorphisms.
   Then the induced map \[F:\operatorname{Cone}_{\bullet}(f)\rightarrow \operatorname{Cone}_
   {\bullet}(\widetilde{f}),\,\,(\alpha,\beta)\mapsto(\Phi(\alpha),\Psi(\beta))\]
   is a quasi-isomorphism.\end{prop}
   \begin{proof}The map $F$ is a chain map
   since,\begin{equation*}\begin{split}
   \partial F(\alpha,\beta)&=\partial(\Phi(\alpha),\Psi(\beta))\\
   & =(\partial\Phi(\alpha),\widetilde{f}(\Phi(\alpha))-\partial\Psi(\beta))\\
   & =(\Phi(\partial(\alpha),\Psi(f(\alpha)-\partial\beta))\\
   & =F(\partial\alpha,f(\alpha)-\partial\beta)\\
   & =F\partial(\alpha,\beta).\end{split}\end{equation*}The chain
   map $F$ fits into a commutative diagram,
   \\*\[\begin{CD}
  0 @>>> Y_{\bullet}  @>>>  \Co_{\bullet}(f) @>>> X_{\bullet-1} @>>> 0  \\
   @.  @VVV  @VVV  @VVV  @.\\
   0 @>>>\widetilde{Y}_{\bullet}@>>>\Co_{\bullet}(\widetilde{f})@>>>\widetilde{X}_{\bullet-1}
   @>>>0
   \end{CD}\]\\*Since $\Phi$ and $\Psi$ are quasi-isomorphisms, so is
   $F$ by Lemma \ref{1}.
   \end{proof}
\begin{prop}\label{3}For any chain map $f_\bullet:X_{\bullet}\rightarrow Y_{\bullet}$ there is a long exact sequence
\[\cdots\rightarrow H_{n-1}(kerf)\overset{j}\rightarrow H_n(f)\overset{k}\rightarrow H_n(coker f)\overset{\delta}
\rightarrow H_{n-2}(kerf)\rightarrow
H_{n-1}(f)\rightarrow\cdots\]where $j$, $k$ and the connecting
homomorphism $\delta$ are defined by\begin{eqnarray*}j[\theta]
&=& [(\theta,0)] \\ k[(\theta,\eta)] &=& [\eta\, mod\, f(X)] \\
\delta[(\eta\, mod\, f(X))] &=& [\partial\theta]\in
H_{n-2}(kerf).\end{eqnarray*} Here, $\eta\in Y_n$ and
$\partial\eta=f(\theta)$ for some $\theta\in X_{n-1}$. In
particular, if $f$ is an injection $H_n(f)=H_n(cokerf)$, and if it
is onto $H_n(f)=H_{n-1}(kerf)$.
\end{prop}\begin{proof} Let $\widetilde{f}_{\bullet}:X_{\bullet}\rightarrow
im(f_\bullet)\subseteq Y_\bullet$ be the chain map $f_\bullet$,
viewed as a map into the sub-complex
$f_\bullet(X_\bullet)\subseteq Y_\bullet$. We have the following
short exact sequence
\[0\rightarrow \operatorname{Cone}_n(\widetilde{f})\overset{i}\rightarrow \operatorname{Cone}_n(f)\overset{k}
\rightarrow coker(f_n)\rightarrow 0\] where $k$ is as above and
$i$ is the inclusion map. Therefore we get a long exact sequence
\begin{equation}
\begin{split}\cdots\rightarrow
H_n(\widetilde{f})\overset{i}\rightarrow
H_n(f)\overset{k}\rightarrow H_n(cokerf) \rightarrow
H_{n-1}(\widetilde{f})\rightarrow\cdots.\end{split}\end{equation}
Let $\widetilde{f}'_\bullet:X_\bullet/kerf_\bullet\rightarrow
im(f_\bullet)$ be the map induced by $f$. Notice that since
$\widetilde{f}'$ is an isomorphism, therefore
$H_\bullet(\widetilde{f}')=0$. By using the long exact sequence
corresponding to the short exact sequence\[0\rightarrow
kerf_{\bullet-1}\overset{\widetilde{j}}\rightarrow
\operatorname{Cone}_\bullet(\widetilde{f})\overset{\pi}
\rightarrow \operatorname{Cone}_\bullet(\widetilde{f}')\rightarrow
0\]where $\widetilde{j}(\theta)=(\theta,0)$, and
$\pi(\theta,\eta)=(\theta\,mod\,kerf,\eta)$, we see that
$\widetilde{j}$ is a quasi-isomorphism. Since
$j=i\circ\widetilde{j}$, we obtain the long exact sequence
\[\cdots\rightarrow H_{n-1}(kerf)\overset{j}\rightarrow H_n(f)\overset{k}\rightarrow
H_n(cokerf)\rightarrow H_{n-2}(kerf)\rightarrow\cdots.\]To find
connecting homomorphism, assume $[\eta$ mod $f(X)]\in H_n(cokerf)$
for $\eta\in Y_n$. Then $\partial\eta\in f(X)$, i.e.,
$\partial\eta=f(\theta)$ for some $\theta\in X_{n-1}$. Since
\[f(\partial\theta)=\partial f(\theta)=\partial
\partial\eta=0\]then $\partial\theta\in ker(f)$. Also
$k(\theta,\eta)= \eta$ mod $f(X)$ and
$j(\theta)=i\circ\widetilde{j}(\partial\theta)=i(\partial\theta,0)=(\partial\theta,0)=\partial(\theta,\eta)$.
Thus, we have  \[\delta[(\eta\, mod\, f(X))]=[\partial\theta]\in
H_{n-2}(kerf).\]
\end{proof}
\goodbreak
\section{Algebraic Mapping Cone for Co-chain Complexes }If $f^\bullet:X^\bullet\rightarrow Y^\bullet$ be a co-chain map
between co-chain complexes, the algebraic mapping cone of $f$ is
defined as a co-chain complex $\operatorname{Cone}^{\bullet}(f)$
where
\[\operatorname{Cone}^n(f)=Y^{n-1}\oplus X^n\]with the differential
\[d(\alpha,\beta)=(f(\beta)-d\alpha
,d\beta)\]Since $d^2=0$, we can consider the cohomology of this
co-chain complex. Define relative cohomology of $f^{\bullet}$ to
be
$$H^n(f):=H^n(\operatorname{Cone}^{\bullet}(f)).$$\begin{rem}
Any cochain complex $(X^{\bullet},d)$ may be viewed as a chain
complex $(\widetilde{X}_{\bullet},\partial)$ where
$\widetilde{X}_n=X^{-n}$ and $\partial_n=d^{-n}\quad(n\in
\mathbb{Z})$. This correspondence takes cochain maps
$f^\bullet:X^{\bullet}\rightarrow Y^{\bullet}$ into chain maps
$\widetilde{f}_\bullet:\widetilde{X}_{\bullet}\rightarrow
\widetilde{Y}_{\bullet}$ where $\widetilde{f}_n=f^{-n}$ and
identifies $\Co(\widetilde{f})$ and $\widetilde{\Co(f)}$ up to a
degree shift:\[\begin{array}{c}
   \widetilde{\Co(f)}_n=\Co(f)^{-n}=Y^{-n-1}\oplus
X^{-n}\\
  \Co(\widetilde{f})_n=\widetilde{X}_{n-1}\oplus\widetilde{Y}_n=X^{-n+1}\oplus
Y^{-n}. \\
\end{array}\]Thus, $\widetilde{\Co(f)}_{n}\cong
\Co(\widetilde{f})_{n+1}.$\newline \indent Using this
corresondence, results for the mapping cone of chain maps directly
carry over to cochain maps.\end{rem}
\goodbreak
\section{Kronecker Pairing}For a chain complex $X_{\bullet}$, the dual
co-chain complex $(X')^{\bullet}$ is defined by
$(X')^{n}=Hom(X_n,R)$ with the dual differential.\begin{prop} Let
$f_\bullet:X_{\bullet}\rightarrow Y_{\bullet}$ be a map between
chain complexes and $(f')^{\bullet}:(Y')^{\bullet}\rightarrow
(X')^{\bullet}$ be its dual cochain map . Then the bilinear
pairing
\[\Co^n(f')\times \Co_n(f)\rightarrow R\]given by the fomula\[\langle(\alpha,
\beta),
(\theta,\eta)\rangle=\langle\alpha,\theta\rangle-\langle\beta,\eta\rangle\]for
$(\alpha,\beta)\in \operatorname{Cone}^n(f')$ and
$(\theta,\eta)\in \operatorname{Cone}_n(f)$, induces a pairing in
cohomology/homology \[H^n(f')\times H_n(f)\rightarrow R.\]
\end{prop}
\begin{proof}
It is enough to show that a cocycle paired with a boundary is zero
and a coboundary paired with a cycle is zero. Let $(\alpha,
\beta)=\partial (\alpha', \beta')$ and $\partial (\theta,\eta)=0$.
Therefore by definition
\[\alpha=f'\beta'-d\alpha',\, \beta=d
\beta'\]and\[\partial \eta=f(\theta), \,\partial \theta=0 .\]
\begin{equation}
\begin{split}
\langle(\alpha,\beta), (\theta,\eta)\rangle &=\langle
\alpha,\theta\rangle-\langle\beta,\eta\rangle\\&=\langle
f'\beta',\theta\rangle-\langle d \alpha',\theta\rangle-\langle d
\beta',\eta\rangle\\&=\langle f'\beta',\theta\rangle-\langle
\alpha',\partial\theta\rangle-\langle
\beta',\partial\eta\rangle\\&=\langle
f'\beta',\theta\rangle-\langle\beta',f(\theta)\rangle\\&=0.
\end{split}
\end{equation} Similarly we can prove that a co-boundary paired with a cycle is
zero.
\end{proof}
\begin{lem}If $f_{\bullet}:X_{\bullet}\rightarrow Y_{\bullet}$ is a chain map
, and $(f')^{\bullet}:(Y')^{\bullet}\rightarrow (X')^{\bullet}$ be
its dual cochain map, then
$\operatorname{Cone}^{\bullet}(f')=(\operatorname{Cone}_{\bullet}(f))'$.\end{lem}\begin{proof}
Notice that
$\operatorname{Cone}^{n}(f')=(\operatorname{Cone}_{n}(f))'=(X^{n-1})'\oplus
(Y^n)'$. It follows from definitions that
\[\langle d(\alpha,\beta),(\theta,\eta)\rangle=
\langle(\alpha,\beta),\partial(\theta,\eta)\rangle.\]Therefore
differential of $\operatorname{Cone}^{n}(f')$ is dual of
differential of $\operatorname{Cone}_{n}(f)$.
\end{proof}
\goodbreak
\section{Singular, De Rham, \v{C}ech Theory}In this Section, we fix two manifolds
$M$ and $N$ and a map $\Phi\in
 C^{\infty}(M,N)$.\\  \textbf{Singular} \textbf{relative}
\textbf{homology}:  Consider the push-forward map
$\Phi_{\ast}:S_q(M, R)\rightarrow S_{q}(N, R)$, where R is a
commutative ring and $S_q(M, R), \, S_q(N, R)$ are the singular
chain complexes of M and N respectively. Singular relative
homology is the homology of the chain complex
$\Co_{\bullet}(\Phi_*)$, and is denoted $H_{\bullet}(\Phi,R)$.
\\ \textbf{Singular} \textbf{relative} \textbf{cohomology}:
Consider the pull-back map $\Phi^{\ast}:S^q(N, R)\rightarrow
S^{q}(M, R)$, where R is a commutative ring and $S^q(M, R), \,
S^q(N, R)$ are the singular co-chain complex of M and N
respectively. Singular relative cohomology is the cohomology of
the co-chain complex $\Co^{\bullet}(\Phi^*)$, and is denoted
$H^{\bullet}(\Phi,R)$.
\\ \textbf{De} \textbf{Rham} \textbf{relative}
\textbf{cohomology}: For $\Phi\in C^{\infty}(M, N)$,  consider the
pull back-map
\[\Phi^*:\Omega^{q}(N)\rightarrow \Omega^q(M)\] between differential co-chain complexes. We denote
the cohomology of $\Co^{\bullet}(\Phi^*)$
 by $H_{dR}^\bullet(\Phi)$ and will call it de Rham relative cohomology
 . \\ \textbf{\v{C}ech} \textbf{relative}
\textbf{cohomology}: Let $A$ be a R-module, and
$\mathcal{U}=\{U_{\alpha}\}$ be a good cover of a manifold $M$,
i.e., all the finite intersections are contractible. For any
collection of indices $\alpha_0,\cdots,\alpha_p$ such that
$U_{\alpha_0}\cap\cdots\cap U_{\alpha_p}\neq\emptyset$, let
\[U_{\alpha_0\cdots\alpha_p}=U_{\alpha_0}\cap\cdots\cap
U_{\alpha_p}.\]A \v{C}ech-p-cochain $f\in
\check{C}^p(\mathcal{U},A)$ is a function
\[f=\coprod_{\alpha_0\cdots \alpha_p} f_{\alpha_0\cdots \alpha_p}: \coprod_{\alpha_0\cdots
\alpha_p}U_{\alpha_0\cdots \alpha_p}\rightarrow A\] where
$f_{\alpha_0\cdots\alpha_p}$ is locally constant and
anti-symmetric in indices. The differential is defined by
\[(d f)_{\alpha_0\cdots\alpha_{p+1}}=\sum_{i=0}^{p+1}(-1)^if_{\alpha_0\cdots\hat{\alpha}_i
\cdots\alpha_{p+1}}\]where the hat means that the index has been
omitted. Since $d\circ d=0$, we can define \v{C}ech cohomology
groups with coefficients in $A$
by\[\check{H}^p(M,A):=H^p(\check{C}(\mathcal{U},A)).\] Let
$\mathcal{U}=\{{U}_{i}\}_{i\in I} , \,
\mathcal{V}=\{{V}_{j}\}_{j\in J}$ be good covers of M and N
respectively such that there exists a map $r:I\rightarrow J$ with
$\Phi({U}_{i})\subseteq {V}_{r(i)}$. Let $\check{C}^\bullet(M, A),
\check{C}^\bullet(N, A)$ be the \v{C}ech complexes for given
covers, where A is an R-module. Using the pull-back map
$\Phi^{\ast}:\check{C}^\bullet(N, A)\rightarrow
\check{C}^\bullet(M,A)$, we define the relative \v{C}ech
cohomology to be the cohomology of $\Co^{\bullet}(\Phi^*)$. Denote
this cohomology by
 $\check{H}^{\bullet}(\Phi, A)$. \newline \indent Suppose that $\underline{A}$
 is one of the sheaves (~\cite{MR95d:14001},~\cite{MR94b:57030})
 $\underline{\mathbb{Z}}, \underline{\mathbb{R}},
 $\underline{U(1)}, $\underline{\Omega^q}$. Denote the space of
k-cochains of the sheaf $\underline{A}$ on $M$ and $N$
respectively by $C^k(M,\underline{A})$ and $C^k(N,\underline{A})$.
Here the differential is defined as above. Again, we have an
induced map
\[ \Phi^*:C^{k}(N,\underline{A})\rightarrow
C^k(M,\underline{A})\]Denote the  cohomology of
$\Co^\bullet(\Phi^*)$ by
$H^*(\Phi,\underline{A})$.\begin{thm}There is a canonical
 isomorphism $H_{dR}^{n}(\Phi)\cong
 H^{n}(\Phi,\mathbb{R})$.\end{thm}
 \begin{proof} Let $S_{sm}^{\bullet}(M,\mathbb{R})$ and $S_{sm}^{\bullet}(N,\mathbb{R})$ be the
 smooth singular cochain complex of $M$
 and $N$ respectively ~\cite{MR83i:57016}. Consider the following
  diagram :  \\*\[\begin{CD}
\Omega^{n}(N)  @>\Phi^*>>\Omega^{n}(M)\\
@Vg^{n}VV   @Vf^{n}VV\\
S_{sm}^{n}(N, \mathbb{R})  @>\Phi^*>>S_{sm}^{n}(M,\mathbb{R})
\end{CD}\] \\*where $f^n$ is  defined by
$f^n(\omega):\sigma\mapsto\int_{\Delta_n}{\sigma^*\omega}$, for
$\omega\in \Omega^n(M)$ and $\sigma\in S_{n}^{sm}(M)$ is a smooth
singular n-simplex. $g^{\bullet}$ is defined in similar fashion.
From these definitions, it is clear that the diagram commutes.
$f^{\bullet}$ and $g^{\bullet}$ are quasi-isomorphisms by de Rham
Theorem ~\cite{MR83i:57016}. Define
$k^{\bullet}:\Omega^{\bullet}(\Phi,\mathbb{R})\rightarrow
S_{sm}^{\bullet}(\Phi,\mathbb{R})$ by $k^{\bullet}(\alpha,
\beta)=(f^{\bullet-1}(\alpha), g^\bullet(\beta))$. We can use
Proposition \ref{2} and deduce that $k^\bullet$ is a
quasi-isomorphism. There is a co-chain map
\[l^{\bullet}:S^{\bullet}(M,\mathbb{R})\rightarrow
S^{\bullet}_{sm}(M,\mathbb{R})\]given by the dual of the inclusion
map in chain level. In (~\cite{MR84k:58001}, p.196) it is shown
that $l^{\bullet}$ is a quasi-isomorphism. Therefore if we use
Proposition \ref{2} again, we get
$$H^{n}(\Phi,\mathbb{R})\cong
H_{sm}^{n}(\Phi,\mathbb{R}).$$Together we have
$H^{\bullet}(\Phi,\mathbb{R})\cong
H_{dR}^{\bullet}(\Phi,\mathbb{R})$.

\end{proof}
 \begin{thm}For $\Phi\in C^{\infty}(M,N)$, there is an
isomorphism $H_{dR}^{q}(\Phi)\cong \check{
H}^{q}(\Phi,\mathbb{R})$.\end{thm}
\begin{proof}Let $\mathcal{U}=\{{U}_{i}\}_{i\in I} , \,
\mathcal{V}=\{{V}_{j}\}_{j\in J}$ be good covers of M and N
together with a map $r:I\rightarrow J$ such that
$\Phi({U}_{i})\subseteq {V}_{r(i)}$. Define the double complex
$E^{p,q}(M)=\check{C}^p(M,\Omega^q)$ where $
\check{C}^p(\mathcal{U},\Omega^q)$ is the set of $q$-forms
$\omega_{\alpha_0\cdots\alpha_p}\in
\Omega^q(U_{\alpha_0\cdots\alpha_p})$ anti-symmetric in indices
with the differential $d$ defined as before. Let
$E^n(M)=\bigoplus_{p+q=n}E^{p,q}(M)$ be the associated total
complex. The map $\Phi:M\rightarrow N$ induces chain maps
$\Phi^*:E^n(N)\rightarrow E^n(M)$. We denote the corresponding
algebraic mapping cone by $E^n(\Phi)$. The inclusion
$\check{C}^n(M,\mathcal{U})\rightarrow E^n(M)$ is a
quasi-isomorphism (~\cite{MR83i:57016}, p. 97). We have a similar
quasi-isomorphism for $N$ and since inclusion maps commute with
pull-back of $\Phi$, we get a quasi-isomorphism
$\check{C}^n(\Phi)\rightarrow E^n(\Phi)$. Thus we get isomorphism
\begin{equation}\label{eq1}\check{H}^n(\Phi,\mathbb{R})\cong
H^n(E(\Phi)).\end{equation} The map $\Omega^n(M)\rightarrow
E^{0,n}(M)\subset E^n(M)$ given by restrictions of forms
$\alpha\mapsto \alpha|_{U_i}$ is a quasi-isomorphism
(~\cite{MR83i:57016}, p. 96). Again these maps commute with pull
back, and hence define a quasi-isomorphism
$\Omega^n(\Phi)\rightarrow E^n(\Phi)$ which means
\begin{equation}\label{eq2}H_{dR}^n(\Phi)\cong
H^n(E(\Phi)).\end{equation} By combining Equation \ref{eq1} and
Equation \ref{eq2}, we get
$\check{H}^{\bullet}(\Phi,\mathbb{R})\cong H_{dR}^{\bullet}(\Phi)$
\end{proof}\begin{rem} A modification of this argument, working
instead with the double complex $C^p(M,S^q)$ given by collection
of $S^q(U_{\alpha_0\cdots\alpha_p})$ gives isomorphism between
\v{C}ech relative cohomology and singular relative cohomology with
integer coefficients, hence $$\check{H}^q(\Phi,\mathbb{Z})\cong
H_S^q(\Phi,\mathbb{Z}).$$
\end{rem}
\goodbreak
\section{Topological Definition of Relative Homology}Let
$\Phi:M\rightarrow N$
 be an inclusion map,
then the push-forward map $\Phi_*:S_\bullet(M,R)\rightarrow
S_\bullet(N,R)$ is injection. Proposition \ref{3} shows that
$H_{\bullet}(\Phi)\cong H_\bullet(S(N)/S(M))=H_{\bullet}(N,M;R)$.
$ H_{\bullet}(N,M;R)$ is known as relative homology. Obviously
this is a special case of what we defined to be a relative
singular homology of an arbitrary map $\Phi:M\rightarrow N$.

Given a continuous map $f:X\rightarrow Y$ of topological spaces,
define mapping cylinder \[\operatorname{Cyl}_f=\frac{(X\times
I)\sqcup Y}{(x, 1)\sim f(x)}\]and mapping cone
\[\Co_f=\frac{\operatorname{Cyl}_f}{X\times \{0\}}
~\cite{MR2002k:55001}.\] Let $\Co(X):= X\times I/X\times
\{0\}.$\clearpage
\begin{figure}\begin{center}\caption{Mapping
Cone}\end{center}\end{figure}\[\includegraphics{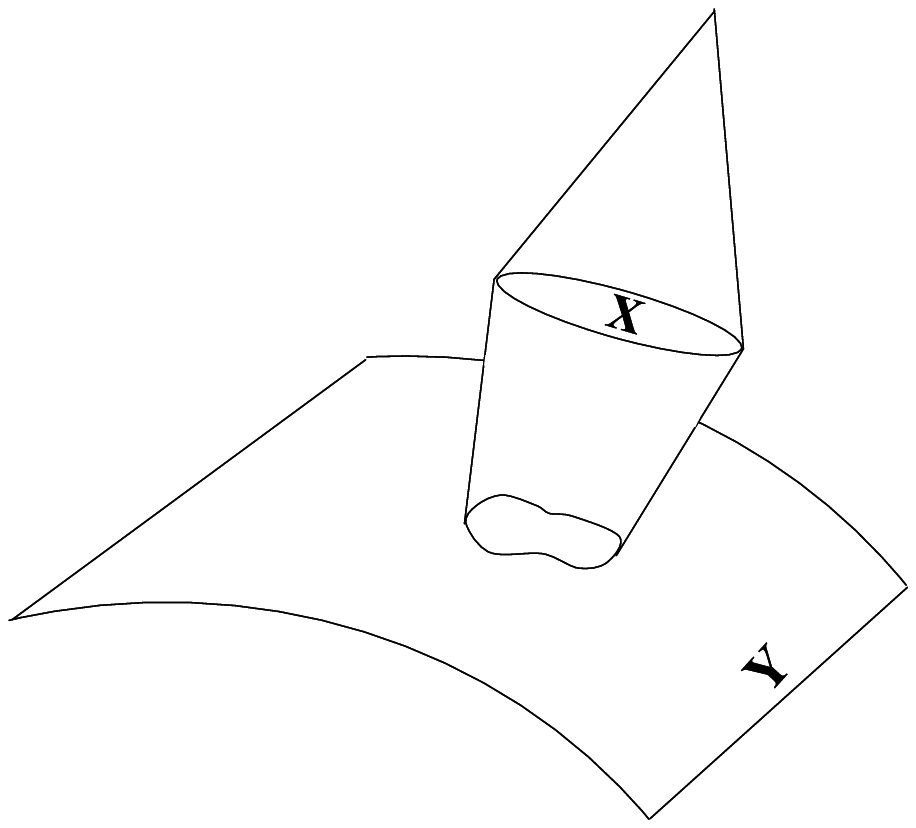}\]There are
natural maps
$$i:Y\hookrightarrow \Co_f,\,\,\,j:\Co(X)\rightarrow
\Co_f.$$Note that $j$ is an inclusion only if $f$ is an inclusion.
There is a canonical map,
$$h: \widetilde{S}_{n-1}(X)\rightarrow \widetilde{S}_{n}(\Co(X))$$
with the property, $h\circ
\partial+\partial\circ h=k$, where $h$ is defined by replacing a singular n-simplex
with its cone and $k:X\hookrightarrow\Co(X)$ is the inclusion map.
Define a map $$l_n:\Co_n(f_*)\rightarrow
S_n(\Co_f),\,\,(x,y)\mapsto j_*(h(x))-i_*(y).$$\begin{thm}
$l_{\bullet}$ is a chain map and is a quasi-isomorphism. Thus,
$$H_n(f)\cong H_n(\Co_f).$$\end{thm}\begin{proof}Recall that
$\partial(x,y)=(\partial x,f_*(x)-\partial y)$. Since
\begin{equation}\begin{split}l(\partial(0,y))+\partial
l(0,y)&= l((0,-\partial y))-\partial i_*y \\
&= i_*(\partial y)-\partial i_*y\\
&=0\end{split}\end{equation}and
\begin{equation}\begin{split}l(\partial(x,0))+\partial
l(x,0)&= l((\partial x,f(x)))+\partial j_*h(x)\\
&= j_*h(\partial x)-i_*f(x)+\partial j_*h(x)\\
&=j_*k_*(x)-i_*f(x)\\
&=0\end{split}\end{equation}therefore $\partial l+l\partial=0$.
Consider diagram \\*\[\begin{CD}
  0 @>>> S_n(\operatorname{Cyl}_f)  @>>>  S_n(\Co_f) @>>> S_n(\Co_f,\operatorname{Cyl}_f) @>>> 0  \\
   @.  @AAA  @A l AA  @AAA  @.\\
   0 @>>>S_n(Y)@>>>\Co_n(f_*)@>>>S_{n-1}(X)
   @>>>0
   \end{CD}\]\\*where the first row corresponds to the pair
   $(\Co_f,\operatorname{Cyl}_f)$ and the right vertical arrow
   comes from$$S_{n-1}(X)\rightarrow S_n(\Co(X),X)\underset{exision}\cong
   S_n(\Co_f,\operatorname{Cyl}_f).$$The diagram commutes, and the
   rows are exact. Since the right and left vertical maps are
   quasi-isomorphisms, hence so is the middle map.\end{proof}
\goodbreak
\section{An Integrality Criterion}
If A and B are R-modules, then any homomorphism
$\kappa:A\rightarrow B $ induces homomorphisms $\kappa:H^n(\Phi,
A)\rightarrow H^n(\Phi, B)$ and $\kappa:H_n(\Phi, A)\rightarrow
H_n(\Phi, B)$. In particular the injection
$\iota:\mathbb{Z}\rightarrow \mathbb{R}$ induces a homomorphism
\[\iota:H^n(\Phi, \mathbb{Z})\rightarrow H^n(\Phi, \mathbb{R})
.\]A class $[\gamma]\in H^n(\Phi, \mathbb{R})$ is called integral
in case $[\gamma]$ lies in the image of the map $\iota$.
\begin{prop}\label{I}A class $[(\alpha, \beta)]\in H^n(\Phi, \mathbb{R})$
is integral if and only if $\int_{\theta}\alpha-\int_{\eta}\beta
\in \mathbb{Z}$ for all cycles
 $(\theta,\eta)\in \Co_n(\Phi, \mathbb{Z})$
. \end{prop}
\begin{proof}
 Consider the following commutative diagram
 \\*\[\begin{CD}0@>>> 0@>>> H^n(\Phi, \mathbb{R})@>{\cong}>>Hom(H_n(\Phi, \mathbb{R}), \mathbb{R})@>>>0\\
@. @AAA @A{\iota}AA @A{\widetilde{\iota}}AA\\
 0@>>>Ext(H_n(\Phi,\mathbb{Z}))@>>>H^n(\Phi, \mathbb{Z})@>{\tau}>>Hom(H_n(\Phi, \mathbb{Z}), \mathbb{Z})@>>>0\\
 \end{CD}\]\\*where
$H^n(\Phi,\mathbb{R})\rightarrow
Hom(H_n(\Phi,\mathbb{R}),\mathbb{R})$ and $\tau$ are pairing given
by integral. The map $\widetilde{\iota}$ is inclusion map,
considering the fact that
$$Hom(H_n(\Phi,\mathbb{R}),\mathbb{R})=
Hom(H_n(\Phi,\mathbb{Z}),\mathbb{R}).$$Thus $[(\alpha, \beta)]\in
H^n(\Phi, \mathbb{R})$ is integral if
$\int_{\theta}\alpha-\int_{\eta}\beta \in \mathbb{Z}$ for all
cycles
 $(\theta,\eta)\in \Co_n(\Phi, \mathbb{Z})$.
 \end{proof}
 \goodbreak
\section{Bohr-Sommerfeld Condition}Let $(N,\omega)$ be a
symplectic manifold. Recall that an immersion $\Phi:M\rightarrow
N$ is isotropic if $\Phi^* \omega=0$. It is called Lagrangian if
furthermore $\dim M=\frac{1}{2}\dim N$. Suppose
$H_1(N,\mathbb{Z})=0$ and $\omega$ is integral. A Lagrangian
immersion $\Phi:M\rightarrow N$ is said to satisfy the
Bohr-Sommerfeld (~\cite{MR85e:58069},~\cite{MR94g:58085})
condition if for all 1-cycles $\gamma\in S_1(M)$
$$\int_D\omega\in \mathbb{Z}\,\,\,\,\,\,\mbox{where}\,\,\,\partial
D=\Phi(\gamma).$$ Note that since $\omega$ is integral, this does
not depend on the choice of $D$. Also, if $\omega=d\theta$ is
exact ( for example for the cotangent bundles ), the condition
means that $$\int_{\gamma}\Phi^*\theta\in
\mathbb{Z}\,\,\,\,\,\,\mbox{for all 1-cycles $\gamma$.}$$ In terms
of relative cohomology, the condition means that $(0,\omega)\in
\Omega^2(\Phi)$ defines an integral class in $H_{dR}^2(\Phi)$. The
interesting feature of this situation is that the forms on $M,N$
are fixed, and defines a condition on the map
$\Phi$.\begin{exmp}Let $N=\mathbb{R}^2,\,M=S^1,\,\omega=dx\wedge
dy,\,\Phi=\mbox{inclusion map}$. Then the immersion $\Phi:
S^1\hookrightarrow \mathbb{R}^2$ satisfies the Bohr-Sommerfeld
condition.\end{exmp}

\chapter{Geometric Interpretation of Integral Relative Cohomology Groups}



Let $\Phi\in C^{\infty}(M,N)$ where M and N are manifolds. Let
$U=\{\emph{U}_{i}\}_{i\in I}$, $V=\{\emph{V}_{j}\}_{j\in J}$ be
good covers of $M$ and $N$ respectively such that there exists a
map $r:I\rightarrow J$ with
$\Phi(\emph{U}_{i})\subseteq\emph{V}_{r(i)}$.
\begin{prop}
$H^q(\Phi, \mathbb{Z})\cong H^{q-1}(\Phi,\underline{U(1)}\,)$ for
$q\geq1$.
\end{prop}
\begin{proof}We have the following long exact sequence
\[\cdots\rightarrow H^{q-1}(M,\underline{\Real})\rightarrow
H^q(\Phi,\underline{\Real})\rightarrow
H^q(N,\underline{\Real})\rightarrow
H^q(M,\underline{\Real})\rightarrow \cdots.\]Since
$H^{\bullet}(M,\underline{\Real})=0$ and
$H^{\bullet}(N,\underline{\Real})=0$, we see that
$H^q(\Phi,\underline{\Real})=0$ for $q>0$. By using the long exact
sequence associated to exponential sequence\[0\rightarrow
\mathbb{Z}\rightarrow
\underline{\mathbb{R}}\overset{\exp}\rightarrow
\underline{U(1)}\rightarrow 0,\,\,\,\,\,\,(\star)\]we can deduce
that $H^q(\Phi,\mathbb{Z})\cong H^{q-1}(\Phi,\underline{U(1)})$
for $q\geq1$.
\end{proof}\goodbreak
\section{Geometric Interpretation of $H^1(\Phi,\mathbb{Z})$}
Let $X$ be a manifold. We say that a function $f\in
C^{\infty}(X,U(1))$ has global logarithm if there exists a
function $k\in C^{\infty}(X,\mathbb{R})$ such that $f=\exp((2\pi
\sqrt{-1}) k)$.

\begin{defn}
Two maps $f,g:X\rightarrow U(1)$ are equivalent if $f/g$ has a
global logarithm.
\end{defn}
The short exact sequence of sheafs $(\star)$ gives an exact
sequence of Abelian groups\[0\rightarrow
H^0(X,\mathbb{Z})\rightarrow
C^{\infty}(X,\mathbb{R})\overset{\exp}\rightarrow
C^{\infty}(X,U(1))\rightarrow H^1(X,\mathbb{Z})\rightarrow 0.\]
This shows there is a one-to-one correspondence between
equivalence classes and elements of $H^1(X,\mathbb{Z})$. We are
looking
 for a geometric realization of $H^1(\Phi,\mathbb{Z})$ for a smooth map $\Phi:M\rightarrow N$. Let\[\mathcal{L}:=
 \{(k,f)|\Phi^*f=\exp((2\pi \sqrt{-1}) k)\}\subset C^{\infty}(M,\mathbb{R})\times
 C^{\infty}(N,U(1)).\]$\mathcal{L}$ has a natural group structure. There is a natural group
homomorphism, \[\tau: C^{\infty}(N,\mathbb{R})\rightarrow
 \mathcal{L}\]where $\tau$
 is defined for $l\in C^{\infty}(N,\mathbb{R})$ by \[\tau(l)=(\Phi^*l,\exp((2\pi \sqrt{-1}) l)).\]
 \begin{defn}We say $(k,f),(k',f')\in \mathcal{L}$ are equivalent if $f/f'=\exp((2\pi \sqrt{-1}) h)$ for
 some function $h\in C^{\infty}(N,\mathbb{R})$ such that\[\Phi^*h=k-k'.\]\end{defn}
 The set of equivalence classes is a group $\mathcal{L}/\tau(C^{\infty}(N,\underline{\Real}))$.\begin{thm}
There exists an exact sequence of groups \[
C^{\infty}(N,\mathbb{R})\overset{\tau}\rightarrow
 \mathcal{L}\rightarrow H^1(\Phi,\mathbb{Z})\rightarrow 0.\]Thus $H^1(\Phi,\mathbb{Z})$
 parameterizes equivalence classes of pairs $(k,f)$.\end{thm}\begin{proof} The first step will be constructing a group homomorphism \[\chi:\mathcal{L}\rightarrow H^1(\Phi,\mathbb{Z}).\]
 Given $(k,f)$, let $l_{j}\in C^{\infty}(V_{j},\mathbb{R})$ be local logarithms
 for $f|_{V_{j}}$, that is $f|_{V_{j}}=\exp((2\pi \sqrt{-1}) l_{j})$. On overlaps,
  $a_{jj'}:=l_{j'}-l_{j} :V_{jj'}\rightarrow \mathbb{Z}$ defines
  a \v{C}ech cocycle in \v{C}$^1(N,\mathbb{Z})$. Let \[b_{i}:= \Phi^*l_{r(i)}-k|_{U_{i}}: U_{i}
  \rightarrow \mathbb{Z}.\]Since $b_i'-b_i=\Phi^*a_{r(i)r(i')}$, so that $(b,a)$ defines a \v{C}ech cocycle in \v{C}$^1(\Phi,\mathbb{Z})$.
  Given another choice of local logarithms $\widetilde{l}_j$, the
  \v{C}ech cocycle changes to
  \[\widetilde{b}_i=b_i+\Phi^*c_{r(i)},\quad
  \widetilde{a}_{jj'}=a_{jj'}+c_{j'}-c_j\]where
  $c_j=\widetilde{l}_j-l_j: V_j\rightarrow \mathbb{Z}$. Thus
  $(\, \widetilde{b},\widetilde{a})=(\,b,a)+ d(0,c)$ and
  $\chi(k,f):=[(b,a)]\in H^1(\Phi,\mathbb{Z})$ is well-defined.
  Similarly, if $(b,a)=d(0,c)$ then the new local logarithms
  $\widetilde{l}_j=l_j-c_j$ satisfy $\widetilde{a}_{jj'}=0$ which
  means that the $\widetilde{l}_j$ patch to a global logarithm
  $\widetilde{l}$. $b_{i}=\Phi^*c_{r(i)}$ implies that
$k|_{U_{i}}=\Phi^*\widetilde{l}_{r(i)}$, which
  means $k=\Phi^*\widetilde{l}$. This shows that the kernel of $\chi$ consists of $(k,f)$ such that
  there exists $l\in C^{\infty}(N,\mathbb{R})$ with $f=\exp((2\pi \sqrt{-1}) l)$ and $k=\Phi^*l$,
  i.e.,
  $ker(\chi)=im(\tau)$.\newline \indent
  Finally, we are going to show that $\chi$ is surjective. Suppose that $(b,a)\in$ \v{C}$^1(\Phi,\mathbb{Z})$ is a cocycle. Then
  \begin{equation}a_{j'j''}-a_{jj''}+a_{jj'}=0\end{equation}\begin{equation}\Phi^*a_{r(i)r(i')}=b_{i'}-b_{i}\end{equation}
  Choose a partion of unity $\sum_{j\in J}h_j=1$ subordinate to
  the open cover $V=\{\emph{V}_{j}\}_{j\in J}$. Define $f_j\in C^{\infty}(V_j,U(1))$
  by\[f_{j}=\exp(2\pi \sqrt{-1}\sum_{p\in J}
  a_{{j}{p}}h_p).\]By applying (2.1.1) on $V_{j}\cap V_{j'}$ we have
  \begin{eqnarray*}f_{j}f_{j'}^{-1}&=&\exp(2\pi \sqrt{-1}\sum_{p\in J}
  a_{{j}{p}}h_p)\exp(-2\pi \sqrt{-1}\sum_{p\in J}
  a_{j' p}h_p)\\
  &=&\exp(2\pi \sqrt{-1}\sum_{p\in J}
  a_{{jj'}}h_p)\\
  &=&1.\end{eqnarray*} Hence the $f_i$ define a map $f\in C^{\infty}(N,U(1))$ such that $f_{|V_{j}}=f_j$.
  Define $k_{i}\in
  C^{\infty}(U_{i},\mathbb{R})$ by\begin{equation}k_{i}=
  \sum_{p\in J}( \Phi^*a_{r(i)p}+b_{i})\Phi^*h_p.\end{equation}Since $b_{i}\in \mathbb{Z}$
  , $\exp((2\pi \sqrt{-1}) k_{i})=\Phi^*f|_{U_{i}}$. We check that on overlaps $U_{i}\cap U_{i'}$,
  $k_{i}-k_{i'}=0$,
 so that the $\{k_{i}\}$ defines a global function $k\in C^{\infty}(M,\Real)$ with $\Phi^*f=\exp((2\pi \sqrt{-1}) k)$. Indeed by applying (2.1.1)
  and (2.1.2) on $U_{i}\cap U_{i'}$ we have
  \[\sum_{p\in J}(\Phi^*a_{r(i) p}+\Phi^*a_{pr(i')}+b_{i}-b_{i'})\Phi^*h_p=\sum_{p\in J}(\Phi^*a_{r(i)r(i')}+b_{i}-b_{i'})\Phi^*h_p=0.\]
  By construction $\chi(b,a)=[(k,f)]$ which shows $\chi$ is surjective.\end{proof}
  \begin{rem}Any $(k,f)\in \mathcal{L}$ defines a $U(1)$-valued
  function on the mapping cone, $\Co_\phi=N\cup_{\Phi}\Co(M)$,
  given by $f$ on $N$ and by $\exp((2\pi \sqrt{-1}) t k)$ on $\Co(M)$. Here
  $t\in I$ is the cone parameter. Hence, one obtains a map $$\mathcal{L}\rightarrow
  H^1(\Co_{\Phi},\mathbb{Z})\cong H^1(\Phi,\mathbb{Z}).$$ This
  gives an alternative way of proving the Theorem
  2.1.1.\end{rem}\goodbreak\goodbreak
\section{Geometric Interpretation of $H^2(\Phi,\mathbb{Z})$}
 Denote the group of Hermitian line bundles over $M$ with $Pic(M)$ and the subgroup of Hermitian line bundles over $M$ which admits
 a unitary section with $Pic_0(M)$.
 Recall that there is an exact sequence of Abelian groups \[0\rightarrow Pic_0(M)\hookrightarrow
 Pic(M)\overset{\delta} \rightarrow H^2(M,\mathbb{Z})\rightarrow 0\]
 defined as follows. For the line bundle $L$ over $M$ with transition maps
 $c_{ii'}\in C^{\infty}(U_{ii'},U(1))$ over
 good cover $\{U_i\}_{i\in I}$ for $M$, $\delta(L)$ is the cohomology class
 of the 2-cocycle $\{a_{ii'i''}: U_{ii'i''}\rightarrow \mathbb{Z}$
 given as
 \[a_{ii'i''}:=\big(\frac{1}{2\pi \sqrt{-1}}(\log c_{i'i''}-\log c_{ii''}+\log c_{ii'})\big)\in\mathbb{Z}.\]
Thus we can say two Hermitian line bundles $L_1$ and $L_2$ over
$M$ are equivalent if and only if $L_1L_2^{-1}$ admits a unitary
section. The exact sequence shows $H^2(M,\mathbb{Z})$
parameterizes the equivalence classes of line
bundles~\cite{MR45:3638}. The class $\delta(L):=c_1(L)$ is called
the first Chern class of $L$. Similarly, for a smooth map
$\Phi:M\rightarrow N$ we are looking for a geometric realization
of $H^2(\Phi,\mathbb{Z})$.

 \begin{defn}
 Suppose $\Phi\in C^{\infty}(M,N)$ and $L_1,L_{2}$ are two Hermitian line bundles
 over N, and $\sigma_{1},\sigma_{2}$ are unitary sections of
 $\Phi^{*}L_1,\Phi^{*}L_{2}$. Then we say
 $(L_1,\sigma_1)$ is equivalent to
 $(L_{2},\sigma_{2})$ if $L_1L_2^{-1}$ admits a unitary section $\tau$
  and there is a map
 $f\in C^{\infty}(M,\mathbb{R})$ such that $(\Phi^{*}\tau)/\sigma_1\sigma_2^{-1}=\exp((2\pi \sqrt{-1}) f)$.
 \end{defn}

  This defines an equivalence relation among $(\sigma, L)$, where
  $L$ is a Hermitian line bundle over N, and $\sigma$ is a
  unitary section of $\Phi^{*}L$. \begin{defn} A relative line bundle for $\Phi\in C^{\infty}(M,N)$ is a
pair $(\sigma, L)$, where $L$ is a Hermitian line bundle over $N$
and $\sigma$ is an unitary section for $\Phi^*L$. Define the group
of relative line bundles \begin{eqnarray*}
  Pic(\Phi)&=&\{(\sigma, L)|L\in Pic(N),
 \sigma \,\mbox{a unitary section of $\Phi^*L$}\}\\
 \mbox{and a
 subgroup of it}\\
Pic_0(\Phi)&=&\{(\sigma, L)\in Pic(\Phi)|\mbox{$\exists$ a unitary
  section}\, \tau\, \mbox{of L and $k$}\in
C^{\infty}(M,\mathbb{R})\\ & &\mbox{ with
}\Phi^*\tau/\sigma=\exp((2\pi \sqrt{-1})
 k)\}.\end{eqnarray*}\end{defn}
 \begin{exmp}Let $(N,\omega)$ be a compact
 symplectic manifold of dimension 2n, and let $L\rightarrow N$ be a line
 bundle with connection $\nabla$ whose curvature is $\omega$,
 i.e., $L$ is a pre-quantum line bundle with connection. A
 Lagrangian submanifold $N$ satisfies the Bohr-Sommerfeld condition if
 there exists a global non-vanishing covariant constant(=flat)
 section $\sigma_M$ of $\Phi^*L$, where $\Phi:M\rightarrow N$ is
 inclusion map (1.6,~\cite{MR94g:58085}). For any Lagrangian submanifold $M$,
 $(\sigma_M, L)\in Pic(\Phi)$.
 \end{exmp}
  \begin{thm}
  There is a short exact sequences of Abelian groups \[0\rightarrow Pic_0(\Phi)
  \rightarrow Pic(\Phi)\rightarrow H^2(\Phi,\mathbb{Z})\rightarrow 0.\]
  Thus $H^2(\Phi,\mathbb{Z})$ parameterizes
  the set of equivalence classes of pairs $(\sigma, L)$.
  \end{thm}

\begin{proof}
We can identify $H^2(\Phi,\mathbb{Z})$ with
$H^1(\,\,\Phi\,,\,\underline{U(1)}\,\,)$ by Proposition 2.0.1. Let
$(\sigma, L)\in Pic(\Phi)$. Let $\{V_{j}\}_{j\in J}$ be a good
cover of N and $\{U_i\}_{i\in I}$ be good cover of $M$ such that
there exist a map $r:I\rightarrow J$ with $\Phi(U_{i})\subseteq
V_{r(i)}$. Choose unitary sections $\sigma_j$ of $L|_{V_{j}}$. The
corresponding transition functions for $L$ are\[g_{jj'}\in
C^{\infty}(V_{jj'},U(1))\quad j,j'\in J,\quad
g_{jj'}\sigma_{j}={\sigma}_{j'}\quad on \quad V_{jj'}\] Define
$f_i=\Phi^*(\sigma_{r(i)})/\sigma$ on $U_{i}$. Then
\begin{equation}\begin{split}f_if_{i'}^{-1}&=(\Phi^*(\sigma_{r(i)})/\sigma).(\Phi^*(\sigma_{r(i')})/\sigma)^{-1}\\&=\Phi^*g_{r(i)r(i')}.\end{split}\end{equation}
Since \[(\delta g)_{r(i)r(i')r(i'')}=1,\]then
$(f_{i},g_{r(i)r(i')})$ is a cocycle in
\v{C}$^1(\Phi,\mathbb{Z})$. If we change local sections
$\sigma_j$, $j\in J$, then $(f_i,g_{r(i)r(i')})$ will shift by a
co-boundary. Define $$\chi:Pic(\Phi)\longrightarrow
H^{1}(\,\,\Phi\,,\,\underline{U(1)}\,\,),\qquad (\sigma,
L)\mapsto[(f,g)].$$\\ \indent To find the kernel of $\chi$,
suppose that $(f,g)=\delta(t,c)$. Thus $g=\delta c$ and
$f=\phi^{*}(c)\exp(2 \pi i h)^{-1}$, where $h$ is the global
logarithm of $t$. Define local section $\tau_j:=\sigma_j/c_j$ on
$V_j$. Since on $V_{jj'}$
\[\sigma_j/c_j=\sigma_{j'}/c_{j'},\] then we obtain a global
section $\tau$. On the other hand,
\[\Phi^*\sigma_{r(i)}/\sigma=f_i=\phi^{*}c_{r(i)}\exp((2\pi \sqrt{-1}) h)^{-1}.\]
Therefore $\phi^*\tau/\sigma=\exp((2\pi \sqrt{-1}) h)^{-1}.$ This
exactly shows that the kernel of $\chi$ is $Pic_0(\Phi)$.\newline
\indent Next we are going to show $\chi$ is onto. Let
$(f_i,g_{jj'})\in C^1(\Phi,\underline{U(1)}\,\,)$ be a cocycle.
Pick a line bundle $L$ over $N$ with $g_{jj'}$ corresponding to
local sections $\sigma_j$. $\Phi^*\sigma_{r(i)}/f_i$ defines local
sections for $\Phi^*L$ over $U_{i}$. On $U_{i}\cap U_{i'}$
\[\Phi^*\sigma_{r(i)}/f_i=\Phi^*\sigma_{r(i')}/f_{i'}\]which defines a global section $\sigma$ for $\Phi^*L$.
By construction, $\chi(\sigma, L)=[(f,g)]$. This shows $\chi$ is
onto.
 \end{proof}
 \begin{rem}A relative line bundle $(L,\sigma)$ for the map
 $\Phi:M\rightarrow N$, defines a line bundle over the mapping
 cone, $\Co_{\Phi}=N\cup_{\Phi}\Co(M)$. This line bundle is given
 by $L$ on $N\subset \Co_{\Phi}$ and by the trivial line bundle
 on $\Co(M)$. The section $\sigma$ is used to glue these two
 bundles. Hence, one obtain a map $$Pic(\Phi)\rightarrow
 H^2(\Co_{\Phi},\mathbb{Z})\cong H^2(\Phi,\mathbb{Z}).$$\end{rem}

 \goodbreak
 \section{Gerbes}Our main references for this Section are
 ~\cite{MR2003f:53086},~\cite{D} and ~\cite{1}.\newline \indent
Let $\mathcal{U}=\{U_i\}_{i\in I}$ be an open cover for a manifold
$M$. It will be convenient to introduce the following notations.
 Suppose there is a collection of line bundles $L_{i^{(0)},...,i^{(n)}}$ on
 $U_{i^{(0)},...,i^{(n)}}$.
 Consider the inclusion maps,\[\delta_k:\emph{U}_{i^{(0)},...,i^{(n+1)}}\rightarrow
 \emph{U}_{i^{(0)},...,\widehat{i^{(k)}},...,i^{(n+1)}}\quad (k=0,\cdots,n+1)\]
 and define Hermitian line bundles $(\delta L)_{i^{(0)},...,i^{(n+1)}}$
 over $\emph{U}_{i^{(0)},...,i^{(n+1)}}$ by\[\delta L:=\underset{k=0}
 {\overset{n+1}{\bigotimes}}
 (\delta_k^{*}L)^{(-1)^k}.\]Notice that $\delta(\delta L)$ is canonically trivial.
 If we have a unitary section $\lambda_{i^{(0)},...,i^{(n)}}$ of $L_{i^{(0)},...,i^{(n)}}$ for each
 $U_{i^{(0)},...,i^{(n)}}\neq\emptyset$, then we can define $\delta\lambda$ in a similar fashion.
 Note that $\delta(\delta\lambda)=1$ as a section of trivial line bundle.
\begin{defn}
A gerbe on a manifold $M$ on an open cover
$\mathcal{U}=\{\emph{U}_{i}\}_{i\in I}$ of $M$, is defined by
Hermitian line bundles $L_{ii'}$ on each $\emph{U}_{ii'}$ such
that $L_{ii'}\cong L_{ii'}^{-1}$ and a unitary section $
\theta_{ii'i''}$ of $\delta L$ on
 $\emph{U}_{ii'i''}$ such
   that $\delta\theta=1$ on $\emph{U}_{ii'i''i'''}$. Denote this data
   by $\mathcal{G}=(\mathcal{U},L,\theta)$.\end{defn} We denote the
   set
   of all gerbes on $M$ on the open cover
$\mathcal{U}=\{\emph{U}_{i}\}_{i\in I}$ by $Ger(M,\mathcal{U})$.
Recall that an open cover $\mathcal{V}=\{V_{j}\}_{j\in J}$ is a
refinement of open cover $\mathcal{U}=\{\emph{U}_{i}\}_{i\in I}$
if there is a map $r:J\rightarrow I$ with $V_j\subset U_i$. In
that case we get a map
\[Ger(M,\mathcal{U})\hookrightarrow Ger(M,\mathcal{V}).\]Define
the group of gerbes on $M$ by\[
Ger(M)=\underset{\longrightarrow}\lim \,Ger(M,\mathcal{U}).\] We
can define the product of two gerbes $\mathcal{G}$ and
$\mathcal{G}'$ to be the gerbe $\mathcal{G}\otimes \mathcal{G}'$
consisting of an open cover of $M$, $\mathcal{V}=\{V_{i}\}_{i\in
I}$ (common refinement of open covers of $\mathcal{G}$ and
$\mathcal{G}'$), line bundles $L_{ii'}\otimes L'_{ii'}$ on
$V_{ii'}$ and unitary sections $\theta_{ii'i''}\otimes
\theta'_{ii'i''}$ of $\delta(L\otimes L')$ on
$V_{ii'i''}$.\newline \indent $\mathcal{G}^{-1}$ the dual of a
gerbe $\mathcal{G}$, is defined by dual bundles $L^{-1}_{ii'}$ on
$U_{ii'}$ and sections $\theta^{-1}$ of $\delta(L^{-1})$ over
$U_{ii'i''}$. Therefore we have a group structure on $Ger(M)$. If
$\Phi:M\rightarrow N$ be a smooth map between two manifolds and
$\mathcal{G}$ be a gerbe on $N$ with open cover
$\mathcal{V}=\{V_{j}\}_{j\in J}$, the pull-back gerbe
$\Phi^*\mathcal{G}$ is simply defined on
$\mathcal{U}=\{U_{i}\}_{i\in I}$ where $\Phi(U_{i})\subset
V_{r(i)}$ for a map $r:I\rightarrow J$, line bundles
$\Phi^*L_{r(i)r(i')}$ on $U_{ii'}$ and  unitary sections $\theta$
of $\delta(\Phi^*L)$ on $U_{ii'i''}$. \goodbreak\begin{defn}
(quasi-line bundle): A quasi-line bundle for the gerbe
$\mathcal{G}$ on a manifold $M$ on the open cover
$\mathcal{U}=\{U_{i}\}_{i\in I}$ is given
by:\begin{enumerate}\item A Hermitian line bundle $E_{i}$ over
each $U_{i}$\item Unitary sections $\psi_{ii'}$ of
\[(\delta E^{-1})_{ii'}\otimes L_{ii'}\]such that
$\delta\psi=\theta$.\end{enumerate} Denote this quasi-line bundle
with $\mathcal{L}=(E,\psi).$
\end{defn} \begin{prop}Any two quasi-line bundle over a given gerbe differ by
a line bundle.\end{prop} \begin{proof} Consider two quasi-line
bundles
 $\mathcal{L}=(E,\psi)$ and $\widetilde{\mathcal{L}}=(\widetilde{E},\widetilde{\psi})$ for the
 gerbe $\g=(\mathcal{U},L,\theta)$.
 $\psi_{ii'}\otimes \widetilde{\psi}_{ii'}^{-1}$ is a unitary section for
 \begin{eqnarray*}E_{i'}\otimes
E_{i}^{-1}\otimes L^{-1}_{ii'}\otimes
\widetilde{E}^{-1}_{i'}\otimes \widetilde{E}_{i}\otimes
L_{ii'}&\cong& E_{i'}\otimes E_{i}^{-1}\otimes
\widetilde{E}^{-1}_{i'}\otimes \widetilde{E}_{i}\\&\cong&
E_{i'}\otimes \widetilde{E}^{-1}_{i'}\otimes E_{i}^{-1}\otimes
\widetilde{E}_{i}.\end{eqnarray*}Therefore $E\otimes
\widetilde{E}^{-1}$ defines a line bundle over $M$.\end{proof}

 Denote the group of all
gerbes on $M$ related to the open cover
$\mathcal{U}=\{\emph{U}_{i}\}_{i\in I}$ which admits a quasi-line
bundle by $Ger_0(M,\mathcal{U})$. Define
\[Ger_0(M)=\underset{\longrightarrow}\lim
\,Ger_0(M,\mathcal{U}).\]

\goodbreak
\begin{prop}\label{D.D}
There exists a short exact sequence of groups \[0\rightarrow
Ger_0(M)\hookrightarrow Ger(M)\overset{\chi}\rightarrow
H^{3}(M,\mathbb{Z})\rightarrow 0.\]
\end{prop}
\begin{proof}
Identify $H^{3}(M,\mathbb{Z})$ with $H^{2}(M,\underline{U(1)})$.
Consider the gerbe $\g$ on $M$. Refine the cover such that all
$L_{ii'}$ admit unitary sections $\sigma_{ii'}$. Define
\[t:=(\delta\sigma)\theta^{-1}.\]Thus, $\delta t=1$ which means $t$
is a cocycle. Define\[ \chi(\mathcal{G}):=[t].\] Different
sections shift the cocycle by $\delta$\v{
C}$^{1}(M,\underline{U(1)})$ which shows $\chi$ is well-defined.
Also $\chi(\mathcal{G}\otimes
\mathcal{G}')=[tt']=\chi(\mathcal{G})\chi(\mathcal{G}')$, which
proves $\chi$ is a group homomorphism. Next, we will show that the
kernel of $\chi$ is $Ger_0(M)$. For $\mathcal{G}\in Ger_0(M)$,
choose a quasi-line bundle $\mathcal{L}=(E,\psi)$. Thus,
$t=\delta(\sigma\psi^{-1})$. Hence, $\chi(\mathcal{G})=[t]=1$.
Conversely, if $[t]=1$, then\[t=\delta t'\] and by defining the
new sections $\sigma'=t'\sigma$ we see that
$\delta\sigma'=t\delta\sigma=\theta$ which shows that
$\mathcal{G}$ admits a quasi-line bundle.\newline \indent Finally,
we show that $\chi$ is onto. If $ \mathcal{U}=\{U_{i}\}_{i\in I}$
be an open cover of $M$ and $t_{ii'i''}$ is a cocycle
\v{C}$^{2}(M,\underline{U(1)})$, then define a gerbe $\mathcal{G}$
on $M$ by trivial line bundle $L_{ii'}$ on $U_{ii'}$ and unitary
sections $\sigma_{ii'}$ on $U_{ii'}$. Define $\theta =
t\delta\sigma $. Since $\delta t=1$, then $\delta\theta=1$. By
construction $\chi(\mathcal{G})=[t]$.

\end{proof}\begin{defn}

Let $\mathcal{G}\in Ger(M)$. $\chi(\mathcal{G})\in
H^2(M,\underline{U(1)})\cong H^3(M,\mathbb{Z})$ is called
Dixmier-Douady class of the gerbe $\mathcal{G}$ which we denote it
by D.D.$\g$.
\end{defn} A gerbe admits a quasi-line bundle if and only if its
Dixmier-Douady class is zero by Proposition \ref{D.D}.
\begin{exmp}\label{Principal bundle}Let $G$ be a Lie group and $1\rightarrow
U(1)\rightarrow \widehat{G} \rightarrow G\rightarrow 1$ a central
extension. Suppose $\pi: P\rightarrow M$ is a principal
$G$-bundle. A lift of $\pi: P\rightarrow M$ is a principal
$\widehat{G}$-bundle $\widehat{\pi}: \widehat{P}\rightarrow M$
together with a map $q: \widehat{P}\rightarrow P$ such that
$\widehat{\pi}=\pi\circ q$ and the following diagram commutes:
\\\[\begin{CD}
   \widehat{G}\times \widehat{P}  @>>>   \widehat{P}\\
   @VVV   @V qVV\\
   G\times P@ >>> P
   \end{CD}\]\\*
Suppose that $\{U_{i}\}_{i\in I}$ be an open cover of $M$ such
that $P|_{U_{i}}:=P_i$ has a lift $\widehat{P}_i$. Define
$\widehat{G}$-equivariant Hermitian line bundles
$$E_i=\widehat{P_i}\times_{U(1)}\mathbb{C}\longrightarrow P\mid_{U_i}.$$
Since $U(1)$ acts by weight 1 on $E_{i}$, it acts by weight 0 on
$E_{i}\otimes E_{i'}^{-1}:=E_{ii'}$ on $U_{ii'}$. Therefore $G$
acts on $E_{ii'}$ and $E_{ii'}/G$ is a well-defined Hermitian line
bundle namely $L_{ii'}$. By construction $\delta L$ is trivial on
$U_{ii'i''}$, therefore we can pick trivial section $\theta$ which
obviously satisfy the relation $\delta\theta=1$. This shows the
obstruction to lifting $P$ to $\widehat{P}$ defines a gerbe
$\g$.\newline \indent If $E_i\rightarrow U_i$ define a quasi-line
bundle $\mathcal{L}$ for $\g$, then the line bundles
$\widetilde{E_i}:=E_i\otimes \pi^*L_i^{-1}$ patch together to a
global $\widehat{G}$-equivariant line bundle
$\widehat{E}\rightarrow P$ and the unit circle bundle defines a
global lift $\widehat{P}\rightarrow P$. Conversely, if $P$ admits
a global lift $\widehat{P}$ and
$\widetilde{P_i}=:\widehat{P}\mid_{U_i}$, then $L_{ii'}$ is
trivial which shows the resulting gerbe is a trivial
one.\end{exmp}
\goodbreak
\begin{exmp}Take $N\subset M$ to be an oriented codimension 3
submanifold of an n-oriented manifold $M$. The tubular
neighborhood $U_0$ of $N$ has the form
$P\times_{SO(3)}\mathbb{R}^3$ where $P\rightarrow N$ is the frame
bundle. Let $U_1=M-N$. Then $U_0\cap U_1\cong
P\times_{SO(3)}(\mathbb{R}^3-0)$. Over $(\mathbb{R}^3-0)\cong
S^2\times (0,\infty)$, we have degree 2 line bundle $E$ which is
$SO(3)$ equivariant. Thus,
$$L_{01}:=P\times_{SO(3)}E$$ is a line bundle over $U_0\cap U_1$
which defines the only transition line bundle. Since there is no
triple intersection, this data defines a gerbe over
$M$.\end{exmp}\goodbreak
  \section{Geometric Interpretation of $H^3(\Phi,\mathbb{Z})$}
\begin{defn} A relative gerbe for $\Phi\in C^{\infty}(M,N)$ is a
pair $(\mathcal{L},\mathcal{G})$, where $\g$ is a gerbe over $N$
and $\mathcal{L}$ is a quasi-line bundle for $\Phi^*\g$.\end{defn}
\textbf{Notation}: Let $\Phi\in C^{\infty}(M,N)$. Then
\begin{eqnarray*}Ger(\Phi)&=&\{(\mathcal{L},\mathcal{G})|(\mathcal{L},\mathcal{G})
 \mbox{is a relative gerbe for}\, \Phi\in C^{\infty}(M,N).\}\\
Ger_0(\Phi)&=&\{(\mathcal{L},\mathcal{G})\in
Ger(\Phi)|\mathcal{G}\,\mbox{admits a quasi-line
bundle}\,\mathcal{L'}\mbox{s.th the line bundle}\,
\mathcal{L}\otimes\Phi^*\mathcal{L'}^{-1}\\& &\mbox{admits a
unitary section}\}\end{eqnarray*} \goodbreak\begin{exmp}Consider a
smooth map $\Phi:M\rightarrow N$ with $\dim M\leq 2$. Let $\g$ be
a gerbe on $N$. Since $\Phi^*\g$ admits a quasi-line bundle say
$\mathcal{L}$, $(\mathcal{L},\g)$ is a relative gerbe.\end{exmp}
\begin{thm}
There exists a short exact sequence of Abelian groups
$$0\rightarrow Ger_0(\Phi)\hookrightarrow Ger(\Phi)\overset{\kappa}\rightarrow
H^{3}(\Phi,\mathbb{Z})\longrightarrow 0$$
\end{thm}
\begin{proof}
We can identify $H^{3}(\Phi,\mathbb{Z})\cong
H^{2}(\Phi,\underline{U(1)})$. Let $\{V_{j}\}_{j\in J}$ be a good
cover of N and $\{U_i\}_{i\in I}$ be good cover of $M$ such that
there exist a map $r:I\rightarrow J$ with $\Phi(U_{i})\subseteq
V_{r(i)}$. Let $(\mathcal{L},\mathcal{G})\in Ger(\Phi)$. Refine
the gerbe $\mathcal{G}=(\mathcal{U},L,\theta)$ sufficiently such
that all $L_{jj'}$ admit unitary sections $\sigma_{jj'}$. Then,
define $t_{jj'j''}\in$ \v{C}$^{2}(N,U(1))$ by
\[t:=(\delta\sigma)(\theta)^{-1}.\]Since $\delta\theta=1$
and $\delta(\delta\sigma)=1$, we have $\delta t=1$. Let
$\mathcal{L}=(E,\psi)$ be a quasi-line bundle for
$\Phi^{*}\mathcal{G}$ with unitary sections $\psi_{ii'}$ for line
bundles $\big((\delta E)_{ii'}\big)^{-1}\otimes
\Phi^{*}L_{r(i)r(i')}$. Define $s_{ii'}\in$\v{C}$^1(M,U(1))$ by
$$s_{ii'}:=(\psi_{ii'})^{-1}\big((\delta\lambda)_{ii'}^{-1}\otimes
\Phi^{*}\sigma_{r(i)r(i')}\big)$$where $\lambda_{i}$ is a unitary
section for $E_{i}$ . Now\begin{equation}\begin{split}
(\Phi^*t^{-1})\delta(\Phi^{*}\sigma)&=\Phi^{*}\theta
\\&=\delta\psi
\\&=(\delta s)^{-1}(\delta\delta\lambda^{-1}\otimes\delta \Phi^{*}\sigma)\\&= (\delta s)^{-1}\delta
\Phi^{*}\sigma.\end{split}\end{equation}This proves that $\delta
s=\Phi^*t$. Define the map
\[\kappa:Ger(\Phi)\rightarrow
H^{2}(\Phi,\underline{U(1)}),
\quad\kappa(\mathcal{L},\mathcal{G})= [(s,t)].\] It is
straightforward to check that this map is well-defined, i.e., it
is independent of the choice of $\sigma_{jj'}$ and $\lambda_i$.
Conversely, given $[(s,t)]\in
H^{2}\big(\Phi,\underline{U(1)}\big)$, we can pick $\mathcal{G}$
such that $\theta=t^{-1}(\delta\sigma)$ and define
$$\psi_{ii'}=s_{ii'}^{-1}\big((\delta\lambda^{-1})_{ii'}\otimes
\Phi^{*}\sigma_{r(i)r(i')}\big).$$ Since $\delta s=\Phi^*t$, then
$\mathcal{L}=(E,\psi)$ defines a quasi-line bundle for
$\Phi^{*}\mathcal{G}$. The construction shows
$\kappa(\mathcal{L},\mathcal{G})=[(s,t)]$. Therefore $\kappa$ is
onto. \\ \indent We now show that $\ker(\kappa)=Ger_0(\Phi)$.
Assume $\kappa(\mathcal{L},\mathcal{G})=[(s,t)]$ is a trivial
class. Therefore there exists $(\rho,\tau)\in$\v{
C}$^{1}\big(\Phi,\underline{U(1)}\big)$ such that
$(s,t)=\delta(\rho,\tau)=(\Phi^{*}\tau(\delta\rho)^{-1},\delta\tau)$.
 $t=\delta\tau$ shows that $\mathcal{G}$ admits a
quasi-line bundle$\mathcal{L'}$. Thus, $\mathcal{L}\otimes
\Phi^*\mathcal{L'}^{-1}$ defines a line bundle over $M$. The first
Chern class of this line bundle is given by the cocycle
$s(\Phi^*\tau)^{-1}$. The condition
$s=(\Phi^*\tau)\delta\rho^{-1}$ shows that this cocycle is exact,
i.e., the line bundle $\mathcal{L}\otimes \Phi^*\mathcal{L'}^{-1}$
admits a unitary section. Thus, $ker(\kappa)\subseteq
Ger_0(\Phi)$. Conversely, if $(\mathcal{L},\mathcal{G})\in
Ger_0(\Phi)$ then the above argument, read in reverse, shows that
$(s,t)$ is exact. Hence, $Ger_0(\Phi)\subseteq ker(\kappa)$.
\end{proof}\goodbreak
\begin{rem}A relative (topological) gerbe $(\mathcal{L},\g)\in Ger(\Phi)$
defines a (topological)gerbe over the mapping cone by ``gluing''
the trivial gerbe over $\Co(M)$ with the gerbe $\g$ over $N\subset
\Co_{\Phi}$. Here, the line bundles $E_i$ which define
$\mathcal{L}$ play the role of transition line bundles. For gluing
of gerbes see ~\cite{S}.\end{rem}\begin{exmp}Let $1\rightarrow
U(1)\rightarrow \widehat{G}\rightarrow G\rightarrow 1$ be a
central extension of a Lie group $G$. Suppose $\Phi\in
C^{\infty}(M,N)$ and $Q\rightarrow N$ is a principal $G$-bundle.
If $P=\Phi^*Q\rightarrow M$ admit a lift $\widehat{P}$, then we
get an element of $H^3(\Phi,\mathbb{Z})$.\end{exmp}
\begin{exmp}Suppose $G$ is a compact Lie group. Recall that the
universal bundle $EG\rightarrow BG$ is a (topological) principal
$G$-bundle with the property that any principal $G$-bundle
$P\rightarrow B$ is obtained as the pull-back by some classifying
map $\Phi:B\rightarrow BG$. While the classifying bundle is
infinite dimensional, it can be written as a limit of finite
dimensional bundles $E_nG\rightarrow B_nG$. For instance, if
$G=U(k)$, one can take $E_nG$ the Stiefel manifold of unitary
$k$-frames over the Grassmanian $Gr_{\mathbb{C}}(k,n)$.
Furthermore, if $B$ is given, any $G$-bundle $P\rightarrow B$ is
given by a classifying map $\Phi:B\rightarrow B_nG$ for some fixed
, sufficiently large $n$ depending only on $\dim
B$~\cite{MR94k:55001}.\\ \indent It maybe shown that
$H^3(BG,\mathbb{Z})$ classifies central extension $1\rightarrow
U(1)\rightarrow\widehat{G}\rightarrow G\rightarrow
1$~\cite{MR97j:58157}. For $n$ sufficiently large,
$H^3(B_nG,\mathbb{Z})=H^3(BG,\mathbb{Z})$. Hence, we find that
$H^3(\Phi,\mathbb{Z})$ classifies pairs
$(\widehat{G},\widehat{P})$, where $\widehat{G}$ is a central
extension of $G$ by $U(1)$ and $\widehat{P}$ is a lift of
$\Phi^*EG$ to $\widehat{G}$.\end {exmp}
\goodbreak
\chapter{Differential Geometry of Relative Gerbes}
\section{Connections on Line Bundles}
 Let $L$ be a Hermitian line bundle with Hermitian connection  $\nabla$ over a manifold $M$.
 In terms
 of local unitary sections $\sigma_i$ of $L\mid_{U_i}$ and the
 corresponding transition maps $$g_{ii'}:U_{ii'}\rightarrow
 U(1),$$connection 1-forms $A_{i}$ on $U_{i}$ are defined by
 $\nabla\sigma_i=(2\pi \sqrt{-1}) A_i \sigma_i$. On $U_{ii'}$, we have
 \[(2\pi \sqrt{-1})(A_{i'}-A_{i})=g_{ii'}^{-1}dg_{ii'}.\]
 Hence, the differentials $dA_i$ agree on overlaps.
 The curvature 2-form $F$ is defined by $F|_{U_{i}}=:dA_{i}$.
 The cohomology class of $F$ is independent of the chosen
connection. The cohomology class of $F$ is the image of the Chern
class $c_1(L)\in H^2(M,\mathbb{Z})$ in $H^2(M,\mathbb{R})$. A
given closed 2-form $F\in\Omega^2(M,\mathbb{R})$ arises as a
curvature of some line bundle with connection if and only if $F$
is integral~\cite{MR94b:57030}.\newline \indent The line bundle
with connection $(L,\nabla)$ is called flat if $F=0$.
 In this case, we define the holonomy of $(L,\nabla)$ as follows:
 We assume
 that the open cover $\{U_i\}_{i\in I}$ is a good cover of $M$.
 Therefore $A_{i}=df_{i}$ on $U_{i}$,
  where $f_{i}:U_i\rightarrow \mathbb{R}$ is a smooth map on $U_{i}$. Then,\[d(2\pi \sqrt{-1}(f_{i'}-f_{i})-\log
  g_{ii'})=0.\]Thus, $$c_{ii'}:=(2\pi \sqrt{-1}(f_{i'}-f_{i})-\log
  g_{ii'})$$ are constants. Since $\log g$ is only defined modulo $2\pi \sqrt{-1}\,\mathbb{Z}$,
  so what we have is a collection of constants
  $\widetilde{c}_{ii'}:=c_{ii'}\mod\mathbb{Z}$. Different choices of $f_i$, shift
  this cocyle with a coboundary. The 1-cocycle $\widetilde{c}_{ii'}$ represents a
\v{C}ech class in \v{H}$^1(M,U(1))$ which is called the
\emph{holonomy} of the flat line bundle $L$ with connection
$\nabla$.\newline \indent
 Let $L\rightarrow M$ be
a line bundle with connection $\nabla$, and $\gamma:S^1\rightarrow
M$ a smooth curve. The holonomy of $\nabla$ around $\gamma$ is
defined as the holonomy of the line bundle $\gamma^*L$ with flat
connection $\gamma^*\nabla$. \section{Connections on Gerbes}
\begin{defn}
Let $\g=(\mathcal{U},L,\theta)$ be a gerbe on a manifold $M$. A
gerbe connection on $\g$ consist of connections $\nabla_{ii'}$ on
line bundles $L_{ii'}$ such that
$\big(\delta\nabla\big)_{ii'i''}\theta_{ii'i''}:=
\big(\nabla_{i'i''}\otimes\nabla_{ii''}^{-1}\otimes\nabla_{ii'}\big)\theta_{ii'i''}=0$
together with 2-forms $\varpi_{i}\in\Omega^{2}(\emph{U}_{i})$ such
that on $U_{ii'}, (\delta \varpi)_{ii'}=F_{ii'}=$ the curvature of
$\nabla_{ii'}$. We denote this connection gerbe by a pair
$(\nabla,\varpi)$. \end{defn}Since $F_{ii'}$ is a closed 2-form,
the de Rham differential $\kappa\mid_{U_i}:=d\varpi_{i}$ defines a
global 3-form $\kappa$ which is called the \emph{curvature of the
gerbe connection}. $[\kappa]\in H^3(M,\mathbb{R})$ is the image of
the Dixmier-Douady class of the gerbe under the induced map by
inclusion $$ \iota: H^3(M,\mathbb{Z})\rightarrow
H^3(M,\mathbb{R}).$$A given closed 3-form
$\kappa\in\Omega^2(M,\mathbb{R})$ arises as a curvature of some
gerbe with connection if and only if i$\kappa$ is
integral~\cite{MR2003f:53086}.
\begin{exmp}Suppose $\pi:P\rightarrow B$ is a principal
$G$-bundle, and $$1\rightarrow U(1)\rightarrow
\widehat{G}\rightarrow G\rightarrow 1$$a central extension. In
Example \ref{Principal bundle}, we described a gerbe $\g$, whose
Dixmier-Douady class is the obstruction to the existence of a lift
$\widehat{\pi}:\widehat{P}\rightarrow B$. We will now explain
(following Brylinski~\cite{MR94b:57030}, see
also~\cite{MR1956150}) how to define a connection on this gerbe.
We will need two ingredients:\newline \indent (i) A principal
connection $\theta\in \Omega^1(P,\mathfrak{g})$,\newline \indent
(ii) A splitting $\tau:P\times_{G}
\widehat{\mathfrak{g}}\rightarrow B\times \mathbb{R}$ of the
sequence of vector bundles
$$0\rightarrow B\times \mathbb{R}\rightarrow P\times_{G}
\widehat{\mathfrak{g}}\rightarrow P\times_{G}
\mathfrak{g}\rightarrow 0$$ associated to the sequence of Lie
algebras $0\rightarrow\mathbb{R}\rightarrow
\widehat{\mathfrak{g}}\rightarrow \mathfrak{g}\rightarrow 0.$ For
a given lift $\widehat{\pi}:\widehat{P}\rightarrow B$, with
corresponding projection $q:\widehat{P}\rightarrow P$, we say that
a principal connection $\widehat{\theta}\in
\Omega^1(\widehat{P},\widehat{\mathfrak{g}})$ lifts $\theta$ if
its image under
$\Omega^1(\widehat{P},\widehat{\mathfrak{g}})\rightarrow
\Omega^1(\widehat{P},\mathfrak{g})$ coincide with $q^*\theta$.
Given such a lift with curvature
$$F^{\widehat{\theta}}=d\widehat{\theta}+\frac{1}{2}[\widehat{\theta},\widehat{\theta}]\in
\Omega^2(\widehat{P},\widehat{\mathfrak{g}})_{\mbox{basic}}=\Omega^2(B,P\times_G
\mathfrak{g}),$$let
$K^{\widehat{\theta}}:=\tau(F^{\widehat{\theta}})\in
\Omega^2(B,\mathbb{R})$ be its ``scalar part''. Any two lifts
$(\widehat{P},\widehat{\theta})$ of $(P,\theta)$ differs by a line
bundle with connection $(L,\nabla^L)$ on $B$. Twisting a given
lift $(\widehat{P},\widehat{\theta})$ by such a line bundle, the
scalar part changes by the curvature of the line bundle
\begin{equation}\label{equation}K^{\widehat{\theta}}+\frac{1}{2\pi
\sqrt{-1}}curv(\nabla^L)~\cite{MR94b:57030}.\end{equation}In
particular, the exact 3-form $dK^{\widehat{\theta}}\in
\Omega^3(B)$ only depend on the choice of splitting and the
connection $\theta$. (It does not depend on choice of lift.) In
general, a global lift $\widehat{P}$ of $P$ does not exist.
However, let us choose local lifts
$(\widehat{P_i},\widehat{\theta}_i)$ of $(P\mid_{U_i},\theta)$.
Denote the scalar part of $F^{\widehat{\theta}_i}$ with
$\varpi_i\in \Omega^2(U_i)$, and let $L_{ii'}\rightarrow U_{ii'}$
be the line bundle with connection $\nabla^{L_{ii'}}$ defined by
two lifts $(\widehat{P_i}\mid_{ U_{ii'}},\widehat{\theta}_i)$ and
$(\widehat{P}_{i'}\mid_{U_{ii'}},\widehat{\theta}_{i'})$. By
Equation
\ref{equation},$$(\delta\varpi)_{ii'}=\frac{1}{2\pi\sqrt{-1}}curv(\nabla^{L_{ii'}}).$$On
the other hand, the connection $\delta\nabla^L$ on $(\delta
L)_{ii'i''}=L_{ii'}L_{ii''}^{-1}L_{ii'}$ is just the trivial
connection on the trivial line bundle. Hence, we have defined a
gerbe connection. \newline \indent A quasi-line bundle $(E,\psi)$
with connection $\nabla^E$ for this gerbe with connection gives
rise to a global lift $(\widehat{P},\widehat{\theta})$ of
$(P,\theta)$, where $\widehat{P}\mid_{U_i}$ is obtained by
twisting $\widehat{P}_i$ by the line bundle with connection
$(E_i,\nabla^{E_i})$. The error 2-form is the scalar part of
$F^{\widehat{\theta}}$.\end{exmp}
\begin{defn}Let $\mathcal{G}$ be a gerbe with connection with a
quasi-line bundle$\mathcal{L}=(E,\psi)$. A connection on a
quasi-line bundle consists of connections $\nabla_i^E$ on line
bundles $E_{i}$ with curvature $F^E_i$ such that
\[(\delta\nabla^E)_{ii'}:=\nabla_{i'}^E\otimes(\nabla_i^E)^{-1}\cong
\nabla_{ii'}.\]Also, the 2-curvatures obey $(\delta
F^E)_{ii'}=F_{ii'}$. We denote this quasi-line bundle with
connection by $(\mathcal{L},\nabla^E)$. Locally defined
 2-forms $\omega\mid_{U_i}=\varpi_{i}-F^E_{i}$ patch together
 to define a global 2-form $\omega$ which is called the error 2-form~\cite{1}.\end{defn}
 \begin{rem}\label{qc}The difference
between two quasi-line bundles with connections is a line bundle
with connection, with the curvature equal to the difference of the
error 2-forms.\end{rem}

 Let $\mathcal{G}=(\mathcal{U},L,\theta)$ be a gerbe with connection on $M$.
 Again, assume that $\mathcal{U}$
 is a good cover. Let
 $t\in$\v{C}$^2\big(M,\underline{U(1)}\big)$ be a representative for the Dixmier-Douady
 class of $\g$. Then,
  we have a collection of 1-forms $A_{ii'}\in \Omega^1(U_{ii'})$ and 2-forms
  $\varpi_{i}\in \Omega^2(U_{i})$ such that \[\kappa|_{U_{i}}=d\varpi_{i}\]\[\delta \varpi=dA\]
  \[(2\pi \sqrt{-1})\delta A=t^{-1}dt.\]If $\kappa=0$, we say the gerbe is flat. In this
  case by using Poincar\'e Lemma,
  $\varpi_i=d\mu_i$ on $U_i$ and on $U_{ii'}$,
  $$(\delta\varpi)_{ii'}=d\delta(\mu)_{ii'}=dA_{ii'}.$$Thus, again by Poincar\'e Lemma
  $$A_{ii'}-(\delta\mu)_{ii'}
  =dh_{ii'}.$$ By using $(2\pi \sqrt{-1})\delta A=t^{-1}dt$, we have
  $$d((2\pi \sqrt{-1})\delta h-\log t)=0.$$
 Therefore, what we have is the collection of constants $c_{ii'i''}
  \in$\v{C}$^2(M,\mathbb{R})$. Since $\log$ is defined
  modulo $2\pi \sqrt{-1}\,\mathbb{Z}$, we define$$\widetilde{c}_{ii'i''}:=c_{ii'i''}\mod \mathbb{Z}.$$
  The 2-cocycle $\widetilde{c}_{ii'i''}$ represents
  a \v{C}ech class in \v{H}$^2(M,U(1))$, which we call it the \emph{holonomy of the flat
  gerbe with
  connection}.
  Let $\sigma:\Sigma \rightarrow M$ be a smooth map, where
$\Sigma$ is a closed
  surface. The holonomy of $\mathcal{G}$ around $\Sigma$ is defined as the holonomy of
  the $\sigma^*\mathcal{G}$ of the flat connection gerbe
  $\sigma^*(\nabla,\varpi)$(~\cite{MR2003f:53086},~\cite{MR1932333}).
  \section{Connections on Relative Gerbes}
  Let $\Phi\in C^{\infty}(M,N)$ and
$U=\{\emph{U}_{i}\}_{i\in I}$, $V=\{\emph{V}_{j}\}_{j\in J}$ are
good covers of $M$ and $N$ respectively such that there exists a
map $r:I\rightarrow J$ with
$\Phi(\emph{U}_{i})\subseteq\emph{V}_{r(i)}$. \begin{defn} A
relative connection on a relative gerbe $(\mathcal{L},\g)$ consist
of gerbe connection $(\nabla,\varpi)$ on $\g$ and a connection
$\nabla^E$ on the quasi-line bundle$\mathcal{L}=(E,\psi)$ for the
$\Phi^*\g$.\end{defn}Consider a relative connection on a relative
gerbe $(\mathcal{L},\g)$. Define the 2-form $\tau$ on $M$ by
$$\tau\mid_{U_i}:=\Phi^*\varpi_{r(i)}-F^E_i.$$ Thus, $(\tau,\kappa)\in \Omega^3(\Phi)$ is
a relative closed 3-form which we call it the \emph{curvature of
the relative connection}.\begin{thm}\label{quantization}A given
closed relative 3-form $(\tau,\kappa)\in \Omega^3(\Phi)$ arises as
a curvature of some relative gerbe with connection if and only if
$(\tau,\kappa)$ is integral.\end{thm}\begin{proof}Let
$(\tau,\kappa)\in \Omega^3(\Phi)$ be an integral relative 3-form.
By Proposition
\ref{I},\begin{equation}\label{This}\int_{\alpha}\kappa-\int_{\beta}\tau\in
\mathbb{Z},\end{equation} where $\alpha\subset N$ is a smooth
3-chain and $\Phi(\beta)=\partial \alpha$, i.e.,
$(\beta,\alpha)\in \Co_3(\Phi,\mathbb{Z})$ is a cycle. If $\alpha$
be a cycle then $(0,\alpha)\in \Co_3(\Phi,\mathbb{Z})$ is a cycle.
In this case equation \ref{This} shows that for all cycles
$\alpha\in S_3(N,\mathbb{Z})$,
$$\int_{\alpha}\kappa\in \mathbb{Z}.$$ Therefore we can pick a
gerbe $\g=(\mathcal{V},L,\theta)$ with connection
$(\nabla,\varpi)$ over $N$ with curvature 3-form $\kappa$. Denote
$\tau_i:=\tau\mid_{U_i}$. Define $F_i^E\in\Omega^2(U_i)$ by
$$F_i^E=\Phi^*\varpi_{r(i)}-\tau_i.$$ Let $(\alpha,\beta)\in\Co_3(\Phi,\mathbb{Z})$
be a cycle. Then
\begin{eqnarray*} \int_{\beta}F_i^E &=&
\int_{\beta}(\Phi^*\varpi_{(r(i)}-\tau_i)\\&=&
\int_{\Phi(\beta)}\varpi-\int_{\beta}\tau\\&=&\int_{\alpha}d\varpi-\int_{\beta}\tau\\
&=&\int_{\alpha}\kappa-\int_{\beta}\tau\in
\mathbb{Z}.\end{eqnarray*}Therefore we can find a line bundle
$E_i$ with connection over $U_i$ whose curvature is equal to
$F_i^E$. Over $U_{ii'}$, the curvature of two line bundles
$\Phi^*L_{ii'}$ and $E_{i'}\otimes E_i^{-1}$ agrees. We can assume
that the open cover $\mathcal{U}=\{U_i\}_{i\in I}$ is a good cover
of $M$. Thus, there is a unitary section $\psi_{ii'}$ for the line
bundle $E_{i}\otimes E_{i'}^{-1}\otimes \Phi^*L_{ii'}$ such that
$\delta\psi=\Phi^*\theta$. Therefore we get a quasi-line bundle
$\mathcal{L}=(E,\psi)$ with connection for $\Phi^*\g$. By
construction the curvature of the relative gerbe $(\lc,\g)$ is
$(\tau,\kappa)$. Conversely, for a given relative gerbe with
connection $(\lc,\g)$ we have $\int_{\beta}F_i^E\in \mathbb{Z}$
where $\beta\subset M$ is a 2-cycle which gives us \ref{This}.
\end{proof} Suppose that $\mathcal{G}$ is a gerbe with a flat
connection gerbe $(\nabla,\varpi)$
   on $N$ and $\mathcal{L}$
  a quasi-line bundle with connection for $\Phi^*\mathcal{G}$.
  Since $\kappa=0$, as it explained in previous Section, we get 2-cocyles
  $\widetilde{c}_{ii'i''}$ which represents a cohomology class in
  \v{H}$^2(M,U(1))$. Since $\Phi^*\mathcal{G}$
  is trivializable, we can find a collection of maps $f_{ii'}$ on $U_{ii'}$
  such that $\delta f=\Phi^*t$ where $j=r(i)$ and $j'=r(i')$. Define $k_{ii'}\in \mathbb{R}$
  by the formula
  \[k_{ii'}=:(2\pi \sqrt{-1})\Phi^*h_{ii'}-\log f_{ii'},\]and \[\widetilde{k}_{ii'}:=
  k_{ii'}\mod\mathbb{Z}.\]Thus, \[\Phi^*\widetilde{c}=\delta \widetilde{k} .\]
  Define the \emph{relative holonomy}
   of the pair $(\,\mathcal{G},\mathcal{L})$ by the relative class
   $[(\widetilde{k},\widetilde{c})]\in H^2(\Phi,U(1))$.

   \begin{defn}
   Let the following diagram be commutative:\\*\[\begin{CD}
   S^1   @>i>>     \Sigma\\
   @VV\psi V   @VV\widetilde{\psi}V\\
   M@>\Phi>>N
   \end{CD}\]\\where $\Sigma$ is a closed surface, i is inclusion map and all other maps are smooth. Suppose
   $\mathcal{G}$ is a gerbe with connection on $N$ and $\Phi^*\mathcal{G}$ admits a quasi-line bundle    $\mathcal{L}$ with
   connection. Clearly $\widetilde{\psi}^*\mathcal{G}$ is a flat gerbe and since $i^*\widetilde
   {\psi}^*\mathcal{G}=\psi^*\Phi^*\mathcal{G}$, $i^*\widetilde{\psi}^*\mathcal{G}$ admits a quasi
   line bundle with connection which is equal to $\psi^*\mathcal{L}$ . We can define the
   holonomy of the relative gerbe around the commutative diagram to be
   the same as holonomy of the pair $(\psi^*\mathcal{L},\widetilde{\psi}^*\mathcal{G})$.
  \end{defn}
\section{Cheeger-Simon Differential Characters}In this section, I
develop a relative version of Cheeger-Simon differential
characters. Denote the smooth singular chain complex on a manifold
$M$ by $S_{\bullet}^{sm}(M)$. Let $Z_{\bullet}^{sm}(M)\subseteq
S_{\bullet}^{sm}(M)$ be the sub-complex of smooth cycles. Recall
that a differential character of degree $k$ on a manifold $M$ is a
homomorphism
$$j:Z_{k-1}^{sm}(M)\rightarrow U(1),$$ such that there is a
closed form $\alpha\in\Omega^k(M)$ with
$$j(\partial x)=\exp\big(2\pi \sqrt{-1}\int_{x}\alpha\big)$$for any
$x\in S_k^{sm}(M)$~\cite{csdc}.

 A connection on a line bundle
defines a differential character of degree 2, with $j$ the
holonomy map. Similarly, a connection on a gerbe defines a
differential character of degree 3. In more details:\newline
\indent Any smooth $k$-chain $x\in S_k^{sm}(M)$ is realized as a
piecewise smooth map $$\varphi_x:K_x\rightarrow M,$$where $K_x$ is
a $k$-dimensional simplicial complex~\cite{MR2002k:55001}. Then,
by definition $$\int_{K_x}\alpha= \int_x\alpha, \qquad \alpha\in
\Omega^k(M).$$ Suppose that $y\in Z_2^{sm}(M)$,
$y=\Sigma\epsilon_i\sigma_i$ with $\epsilon_i=\pm 1$. Assume $\g$
is a gerbe with connection over $M$. Since
$H^3(K_y,\mathbb{Z})=0,$ $\varphi_y^*\g$ admits a piecewise smooth
quasi-line bundle $\mathcal{L}$ with connection. That is, a
quasi-line bundles $\mathcal{L}_i$ for all
$\varphi^*\g\mid_{\Delta^k_{\sigma_i}}$, such that they agree on
the matching boundary faces. Let $\omega\in \Omega^2(K_y)$ be the
error 2-form and define
$$j(y):=\exp\big(2\pi \sqrt{-1}\int_{K_{y}}\omega\big).$$Any two
quasi-line bundles differ by a line bundle and hence different
choices for $\mathcal{L}$, changes $\omega$ by an integral 2-form.
This shows that $j$ is well-defined. Assume that $y=\partial x$.
Since the components of $K_x$ with empty boundary will not
contribute, we can assume that each component of $K_x$ has
non-empty boundary. Since $H^3(K_x,\mathbb{Z})=0$, we can choose a
quasi-line bundle with connection for $\varphi^*_x\g$ with error
2-form $\omega$. Let $k$ be the curvature of $\g$, since
$\varphi^*_xk=d\omega$, by stokes' Theorem we have
$$\int_{K_x}k=\int_{K_x}d\omega=\int_{\partial
K_x}\omega=\int_{K_y}\omega.$$This shows that $j$ is a
differential character of degree 3.

\begin{defn}Let $\Phi\in C^{\infty}(M,N)$ be a smooth map between
manifolds. A relative differential character of degree $k$ for the
map $\Phi$ is a homomorphism
$$j:Z_{k-1}^{sm}(\Phi)\rightarrow U(1),$$ such that there is a
closed relative form $(\beta,\alpha)\in\Omega^k(\Phi)$ with
$$j(\partial(y,x))=\exp\big(2\pi \sqrt{-1}(\int_{y}\beta-\int_{x}\alpha)\big)$$for any
$(y,x)\in S_k^{sm}(\Phi)$.\end{defn}\begin{thm}A relative
connection on a relative gerbe defines a relative differential
character of degree 3.\end{thm} \begin{proof}Let $\Phi\in
C^{\infty}(M,N)$ be a smooth map between manifolds and consider a
relative gerbe $(\mathcal{L},\g)$ with connection. Let $(y,x)\in
S^{sm}_1(\Phi)$ be a smooth relative singular cycle, i.e.,
$$\partial y=0$$and $$\Phi_*(y)=\partial x.$$ Let $K_y$ and $K_x$ be
the corresponding simplicial complex, and $$\Phi:K_y\rightarrow
K_x$$ be the induced map. Given a relative connection, choose a
quasi-line bundle $\mathcal{L'}$ for $\varphi_x^*\g$ and a unitary
section $\sigma$ of the line bundle
$H:=\varphi_y^*\mathcal{L}\otimes(\Phi^*\mathcal{L}')^{-1}$. Let
$\widetilde{\omega}\in \Omega^2(N)$ be the error 2-form for
$\mathcal{L}'$ and $A\in \Omega^1(M)$ be the connection 1-form for
$H$ with respect to $\sigma$. Define a map $j$ by
$$j(y,x):=\exp\big(2\pi \sqrt{-1}\big(\int_{K_x}\widetilde{\omega}-\int_{K_y}A\big)\big).$$
If we choose another quasi-line bundle for $\varphi_x^*\g$, the
difference of error 2-forms will be an integral 2-form and
changing the section $\sigma$ will shift connection 1-form $A$ to
$A+A'$, where $A'$ is an integral 1-form. This proves that
$$j:Z_2^{sm}(\Phi)\rightarrow U(1)$$ is well-defined. Let
$k$ be the curvature 3-form for $\g$ and $\omega$ be the error
2-form for $\mathcal{L}$. Then $(\omega,k)\in \Omega^3(\Phi)$ and
\begin{eqnarray*}j\big(\partial(y,x)\big)&=&j(\partial
y,\Phi_*(y)-\partial x)\\
&=&\exp\big(2\pi \sqrt{-1}\big(\int_{K_{(\Phi_*(y)-\partial
x)}}\widetilde{\omega}-\int_{K_{\partial
y}}A\big)\big)\\&=&\exp\big(2\pi
\sqrt{-1}\big(\int_{K_{\Phi_*(y)}}\widetilde{\omega}
-\int_{K_{\partial x}}\widetilde{\omega}-\int_{K_{\partial
y}}A\big)\big)\\&=&\exp\big(2\pi
\sqrt{-1}\big(\int_{K_{y}}(\Phi^*\widetilde{\omega}-dA)
-\int_{K_x}d\widetilde{\omega}\big)\big)\\
&=&\exp\big(2\pi \sqrt{-1}\big(\int_{K_{y}}\omega
-\int_{K_x}k\big)\big).\end{eqnarray*}Thus, $j$ is a relative
differential character in degree 3.\end{proof}
\section{Relative Deligne Cohomology}
Suppose we have a co-chain complex of sheaves
$\underline{A}^{\bullet}$ over a manifold $M$. For an open cover
$\mathcal{U}=\{U_i\}_{i\in I}$, define
A$^{p,q}_{\mathcal{U}}:=$\v{C}$^p(\mathcal{U},\underline{A}^q)$ to
be the \v{C}ech p-cochains with values in the sheaf
$\underline{A}^q$. On A$^{p,q}_{\mathcal{U}}$ we have two
differentials, $d^q:$A$^{p,q}_{\mathcal{U}}\rightarrow$
A$^{p,q+1}_{\mathcal{U}}$ which is induced from the differential
of co-chain complex of sheaves $\underline{A}^{\bullet}$ and
\v{C}ech coboundary map
$\delta:$A$^{p,q}_{\mathcal{U}}\rightarrow$
A$^{p+1,q}_{\mathcal{U}}$. Since $\delta d=d\delta$, $d^2=0$ and
$\delta^2=0$, A$^{\bullet,\bullet}_{\mathcal{U}}$ is a double
complex. We denote the total complex of this double complex by
$T^{\bullet}_{(\mathcal{U},A)}$. A refinement
$\mathcal{V}<\mathcal{U}$ will induce mappings
$T^{\bullet}_{(\mathcal{U},A)}\rightarrow
T^{\bullet}_{(\mathcal{V},A)}$. The \emph{sheaf  hypercohomology
groups} (~\cite{MR94b:57030},~\cite{MR95d:14001}) of $M$ with
values in the cochain complex of sheaves $\underline{A}^{\bullet}$
are defined by
\[H^{\bullet}(M,\underline{A}):=\underset{\rightarrow}\lim H^{\bullet}(T^{\bullet}_
{(\mathcal{U},A)}).\]For any $f\in C^{\infty}(U_i,U(1))$
 define $$d\log
 f:=\frac{1}{2\pi\sqrt{-1}}f^{-1}df\in\Omega^1(U_i).$$
It is known that classes in the first Deligne hypercohomology
group
$$H^1(M;\underline{U(1})\overset{dlog}\longrightarrow
\underline{\Omega}^1)$$ are in a one-to-one correspondence
 with isomorphism classes of line bundle with connection on $M$
 and classes in the second Deligne
hypercohomology group
$$H^2(M;\underline{U(1})\overset{dlog}\longrightarrow
\underline{\Omega}^1\overset{d}\longrightarrow
\underline{\Omega}^2)$$ are in a one-to-one correspondence
 with isomorphism classes of gerbes with connection on
 $M$~\cite{MR94b:57030}.
 Suppose $M$ and $N$ are two manifolds. Fix a map $\Phi\in
 C^{\infty}(M,N)$. Let
$\mathcal{U}=\{{U}_{i}\}_{i\in I} , \,
\mathcal{V}=\{{V}_{j}\}_{j\in J}$ be good covers of M and N
respectively such that there exists a map $r:I\rightarrow J$ with
$\Phi({U}_{i})\subseteq {V}_{r(i)}$. Also assume
$\underline{A}^{\bullet}$ and $\underline{B}^{\bullet}$ are two
cochain complexes of sheaves over $M$ and $N$ respectively and one
is given a homomorphism of sheaf complexes
$\Phi^*:\underline{B}^{\bullet}\rightarrow
\underline{A}^{\bullet}$. Then $\Phi$ induces mappings
map\[(\Phi^*)^{p,q}:B^{p,q}_{\mathcal{V}}\rightarrow
A^{p,q}_{\mathcal{U}}.\]$\Phi^*$ will induce mappings
\[(\widetilde{\Phi^*})^{\bullet}:T^{\bullet}_{(\mathcal{V},B)}\rightarrow
T^{\bullet}_{(\mathcal{U},A)}.\]Denote the algebraic mapping cone
of $(\widetilde{\Phi^*})^{\bullet}$ by
$\Co^{\bullet}_{(\mathcal{U},\mathcal{V})}(\widetilde{\Phi^*})$
and its corresponding cohomology by $H^{\bullet}_{(\m,\n)}(\Phi)$.
If $(\m',\n')$ is a double-refinement of $(\m,\n)$, i.e., $\m'$ is
a refinement of $\m$, $\n'$ is a refinement of $\n$ and there are
maps $r:I\rightarrow J$ with $\Phi({U}_{i})\subseteq
{V}_{r(i)},\quad r':I'\rightarrow J'$ with
$\Phi({U'}_{i})\subseteq {V'}_{r(i)}$, then we will have induced
mappings
\[\Co^{\bullet}_{(\mathcal{U},\mathcal{V})}(\widetilde{\Phi^*})\rightarrow\Co^{\bullet}
_{(\mathcal{U'},\mathcal{V'})}(\widetilde{\Phi^*}).\]Define the
\emph{relative sheaf hypercohomology groups}
$H^{\bullet}(\Phi,\underline{A},\underline{B})$ by
\[H^{\bullet}(\Phi,\underline{A},\underline{B}):=\underset{\underset{(\m,\n)}\rightarrow}\lim H_{(\m,\n)}^{\bullet}(\Phi).\]
\begin{defn}We define relative Deligne cohomology of $\Phi$ by
$$H^{\bullet}(\Phi;\underline{U(1})\overset{dlog}\longrightarrow
\underline{\Omega}^1;\underline{U(1})\overset{dlog}\longrightarrow
\underline{\Omega}^1\overset{d}\longrightarrow
\underline{\Omega}^2)$$ and we will denote it by
$H_{D}^{\bullet}(\Phi)$.\end{defn}
\begin{thm}There is a one-to-one correspondence between classes of
 $H_{D}^{2}(\Phi)$ and isomorphism classes of relative gerbe with
connection.\end{thm}\begin{proof} Let
$\mathcal{G}=(\mathcal{V},L,\theta)$ be a gerbe with connection on
$N$. Assume that $\mathcal{V}$
 is a good cover. Let
 $t\in$\v{C}$^2\big(N,\underline{U(1)}\big)$ be a representative for the Dixmier-Douady
 class of $\g$. Then,
  we have a collection of 1-forms $A_{jj'}\in \Omega^1(V_{jj'})$, 2-forms
  $\varpi_{j}\in \Omega^2(V_{j})$ and 3-form $\kappa
  \in \Omega^3(N)$ such that \[\kappa|_{V_{j}}=d\varpi_{j}\]\[\delta \varpi=dA\]
  \[(2\pi \sqrt{-1})\delta A=t^{-1}dt.\]This defines a class $(\varpi,A,t)\in
H^2(N;\underline{U(1})\overset{dlog}\longrightarrow
\underline{\Omega}^1\overset{d}\longrightarrow
\underline{\Omega}^2).$ Let $\lc$ be a quasi-line bundle with
connection for $\Phi^*\g$. We can use $\Psi$ to get $h\in$
\v{C}$^1\big(M,\underline{U(1)}\big)$ such that \[\Phi^*t=\delta
h.\]Also by using the identity $(\delta
F^E)_{ii'}=(\Phi^*F)_{r(i)r(i')}$, we can find $l\in
\Omega^1(U_{ii'})$ such that
\[\Phi^*A=\delta l\]Therefore,
a relative gerbe $(\lc,\mathcal{G})$ gives a class
$[((l,h),(\varpi,A,t))]\in H^2_D(\Phi)$.\end{proof}
\goodbreak
\section{Transgression}
For a manifold $M$, we denote its loop space by $LM$. In this
Section first we will construct a line bundle with connection over
$LM$ by transgressing a gerbe with connection over $M$. A map
$\Phi\in C^{\infty}(M,N)$ induces a map $L\Phi\in
C^{\infty}(LM,LN)$. We will prove that a relative gerbe with
connection on $\Phi$ produce a relative line bundle with
connection on $L\Phi$ by transgression. In cohomology language it
means that there is a map\[T:H_D^{\bullet}(\Phi)\rightarrow
H_D^{\bullet-1}(L\Phi).\] \begin{prop}(Parallel transportation)
Suppose that $\g$ is a gerbe with connection on $M\times [0,1]$
and $\g_0=\g|_{(M\times \{0\})}$. There is a natural quasi-line
bundle with connection for the gerbe $\pi^*\g_0\otimes\g^{-1}$,
where $\pi$ is the projection map
\[\pi:M\times [0,1]\rightarrow M\times \{0\}.\]\end{prop}\begin{proof}It is obvious that we
can get a quasi-line bundle with connection for the gerbe
$\pi^*\g_0\otimes\g^{-1}$. We will specify a quasi-line bundle
${\lc}_{\g}$ for the gerbe $\pi^*\g_0\otimes\g^{-1}$ by the
following requirements:\newline \indent 1.The pull-back
$\iota^*{\lc}_{\g}$ is trivial, while $\iota$ is inclusion map
\[\iota: M\times \{0\}\hookrightarrow M\times [0,1].\]2. Let
$\eta\in\Omega^3(M\times [0,1])$ be the connection 3-form for
$\pi^*\g_0\otimes\g^{-1}$. We have $\iota^*\eta=0$. Let
$\chi\in\Omega^2(M\times [0,1])$ be canonical primitive of $\eta$
given by transgression. Then we choose the connection on
${\lc}_{\g}$ such that its error 2-form is $\chi$. Any two such
quasi-line bundles differ by a flat line bundle over $M\times
[0,1]$ and which is trivial over $M\times \{0\}$.\end{proof}
\begin{thm}\label{10} A gerbe $\g$ with connection
on $M\times S^1$, induces a line bundle $E_{\g}$ with connection
on $M$. Also, a quasi-line bundle with connection for $\g$ induces
a unitary section of $E_{\g}$.\end{thm}\begin{proof}$M\times
S^1=M\times [0,1]/\sim$, where the equivalence relation is defined
by $(m,0)\sim (m,1)$ for $m\in M$. Therefore
$\pi^*\g_0\otimes\g^{-1}|_{M\times \{1\}/\sim}$ is a trivial gerbe
and $\lc_{\g}|_{M\times \{1\}/\sim}$ is a quasi-line bundle with
connection for this trivial gerbe, i.e., a line bundle with
connection $E_{\g}$ for $M$. If we change $\lc_{\g}$ to another
natural quasi-line bundle with connection, the difference between
two line bundles over $M\times S^1$ is a trivial line bundle. This
shows the assignment $\g \rightarrow E_{\g}$ is well-defined.
\newline \indent If the gerbe $\g$ admits a quasi-line bundle $\lc$, then
$(\pi^*\lc_0)\otimes(\lc^{-1})$ and $\lc_{\g}$ are two quasi-line
bundles for the gerbe $\pi^*\g_{0}\otimes\g^{-1}$, where
$\lc_0=\lc|_{M\times\{0\}}$. Thus
$\pi^*\lc_{0}\otimes\lc^{-1}\otimes(\lc_{\g})^{-1}$ defines a line
bundle over $M\times S^1=M\times [0,1]/\sim$. This line bundle
over $M$ defines a map $s:M\rightarrow U(1)$.
$(\pi^*\lc_0)\otimes(\lc^{-1})|_{M\times\{0\}/\sim}$ is the
trivial line bundle $E$. Since $E_{\g}\otimes E^{-1}=s$, therefore
$E_{\g}$ admits a unitary section.\end{proof}
\begin{rem}Let $X=LM$ and $\g$ is a gerbe with connection $M$. Consider the
evaluation map \[e:LM\times S^1\rightarrow M.\]Then $e^*\g$ will
induce a line bundle with connection on $LM$.\end{rem}\begin{thm}
For a given map $\Phi\in C^{\infty}(M,N)$, a relative gerbe with
connection $\g_{\Phi}$ will induce a relative line bundle with
connection $E_{L\Phi}$.\end{thm}\begin{proof}The relative gerbe
$\g_{\Phi}$ is a gerbe $\g$ on $N$ together with a quasi-line
bundle with connection $\lc$ for the pull-back gerbe $\Phi^*\g$.
The gerbe $\g$ induces a line bundle with connection $E_{\g}$ and
the quasi-line bundle with connection $\lc$ for $\Phi^*\g$ induces
a unitary section $s$ for the line bundle with connection
$(L\Phi)^*E_{\g}$ by Theorem \ref{10}. Thus, the pair $(s,E_{\g})$
defines a relative line bundle with connection
$E_{L\Phi}$.\end{proof}

\def\baselinestretch{1.8}

\def\baselinestretch{1.8}
\goodbreak
\chapter{Pre-quantization of Group-Valued Moment Maps}
\section{Gerbes over a Compact Lie Group }It is
well-known fact that for a compact, simple, simply connected Lie
group the integral cohomology $H^{\bullet}(G,\mathbb{Z})$ is
trivial in degree less than three, while $H^3(G,\mathbb{Z})$ is
canonically isomorphic to $\mathbb{Z}$. The gerbe corresponding to
the generator of $H^3(G,\mathbb{Z})$ is called the basic gerbe
over $G$. In this section, I give an explicit construction of the
basic gerbe for $G=SU(n)$, and of suitable multiples of the basic
gerbe for the other Lie groups. This gerbe plays an important role
in pre-quantization of the quasi-Hamiltonian
$G$-spaces.\subsection{Some Notations from Lie Groups}Let $G$ be a
compact, simple simply connected Lie group and with a maximal
torus $T$. Let $\mathfrak{g}$ and $\mathfrak{t}$ denote the Lie
algebras of $G$ and $T$ respectively. Denote by $\Lambda\subset
\mathfrak{t}$ the integral lattice, given as the kernel of
$$\exp:\mathfrak{t}\rightarrow T.$$
Let $\Lambda^*=Hom(\Lambda,\mathbb{Z})\subset \mathfrak{t}^*$ be
its dual weight lattice. Recall that any $\mu\in \Lambda^*$
defines a homomorphism
$$h_{\mu}:T\rightarrow U(1),\,\exp\xi\mapsto e^{2\pi\sqrt{-1}
\langle\mu,\xi\rangle}.$$ This identifies $\Lambda^*=Hom(T,U(1))$.
Let $\mathcal{R}\subset \Lambda^*$ be the set of roots, i.e., the
non-zero weights for the adjoint representation. Define
$$\mathfrak{t}^{reg}:=\mathfrak{t} \setminus\bigcup_{\alpha\in \mathcal{R}}
ker\alpha.$$ The closures of the connected components of
$\mathfrak{t}^{reg}$ are called Weyl chambers. Fix a Weyl chamber
$\mathfrak{t}_+$. Let $\mathcal{R}_+\subset \Lambda^*$ be the set
of the positive roots, i.e., roots that are non-negative on
$\mathfrak{t}_+$. Then $\mathcal{R}=\mathcal{R}_+\cup
-\mathcal{R}_+$. A positive root is called simple, if it cannot be
written as the sum of positive roots. We denote the set of simple
roots by $\mathcal{S}$.  The set of simple roots
$\mathcal{S}\subset \mathcal{R}_+$ forms a basis of
$\mathfrak{t}$, and
$$\mathfrak{t}_{+}=\{\xi\in\mathfrak{t}\mid\langle\alpha,\xi\rangle\geq 0,
\forall\alpha\in S\}.$$ Any root $\alpha\in \mathcal{R}$ can be
uniquely written as
$$\alpha=\Sigma k_i\alpha_i,\quad k_i\in \mathbb{Z},\alpha_i\in\mathcal{S}.$$
The hight of $\alpha$ is defined by
$\operatorname{ht}(\alpha)=\Sigma k_i$. Since $\mathfrak{g}$ is
simple, there is a unique root $\alpha_0$ with
$\operatorname{ht}(\alpha)\geq \operatorname{ht}(\alpha_0)$ for
all $\alpha\in \mathcal{R}$ which is called the lowest root. The
fundamental alcove is defined as
$$\mathfrak{A}=\{\xi\in\mathfrak{t}_{+}\mid\langle\alpha_0,\xi\rangle\geq-1\}.$$
 The basic inner product
on $\mathfrak{g}$ is the unique invariant inner product such that
$\alpha.\alpha=2$ for all long roots $\alpha$ which we use it to
identify $\mathfrak{g}^*\cong\mathfrak{g}$. The mapping
$\xi\rightarrow \operatorname{Ad} G(\exp\xi)$ is a homeomorphism
from $\mathfrak{A}$ onto $G/\operatorname{Ad} G$, the space of the
cojugacy classes in $G$. Therefore the fundamental alcove
parameterizes conjugacy classes in $G$~\cite{DK}. We will denote
the quotient map by $q:G\rightarrow \mathfrak{A}$.\newline \indent
Let $\theta^L,\theta^R\in \Omega^1(G,\mathfrak{g})$ be the left
and right invariant Maurer Cartan forms. If $L_g$ and $R_g$ denote
left and right multiplication by $g\in G$, then the values of
$\theta^L_{g}$ and $\theta^R_{g}$ at $g$ are given by
\[\theta^L_{g}=dL_{g^{-1}} : TG_g\rightarrow
TG_e\cong\mathfrak{g},\quad \theta^R_{g}=dR_{g^{-1}} :
TG_g\rightarrow TG_e\cong\mathfrak{g}.\]For any $g\in G$,
\[\theta^L_{g}=Ad_g(\theta^R_{g}).\]For any invariant inner product
$B$ on $\mathfrak{g}$, the form
\begin{equation}\label{3-form}\eta:=\frac{1}{12} B(\theta^L,
[\theta^L,\theta^L])\in \Omega^3(G)\end{equation} is bi-invariant
since the inner product is invariant. Any bi-invariant form on a
Lie group is closed, therefore $\eta$ is a closed 3-form. Its
cohomology class represents the generator of
$H^3(G,\mathbb{R})=\mathbb{R}$ if we assume that $G$ is compact
and simple. If in addition $G$ is simply connected, then
$H^3(G,\mathbb{Z})=\mathbb{Z}$ and one can normalize the inner
product such that $[\eta]$ represents an integral generator
~\cite{MR2001i:53140}. We will say that $B$ is the inner product
at level $\lambda>0$ if $B(\xi,\xi)=2\lambda$ for all short
lattice vectors $\xi\in \Lambda$. The inner product at level
$\lambda=1$ is called the basic inner product. (It is related to
the Killing form by a factor $2c_{\mathfrak{g}}$, where
$c_{\mathfrak{g}}$ is the dual Coxeter number of $\mathfrak{g}$.)
Suppose $G$ is simply connected and simple. It is known that the
3-form defined by Equation \ref{3-form} is integral if and only if
its level $\lambda$ of $B$ is an integer.
\subsection{Standard Open Cover of G} Let $\mu_0,\cdots,\mu_d$ be
the vertices of $\mathfrak{A}$, with $\mu_0=0$. Let
$\mathfrak{A}_j$ be the complement of the closed face opposite to
the vertex $\mu_j$. The standard open cover of $G$ is defined by
the pre-images $V_j=q^{-1}(\mathfrak{A}_j)$. Denote the
centralizer of $\exp\mu_j$ by $G_j$. Then the flow-out
$S_j=G_j.\exp(\mathfrak{A}_j)$ is an open subset of $G_j$, and is
a slice for the conjugation action of $G$. Therefore
$$G\times_{G_{j}}S_j=V_j.$$More generally let
$\mathfrak{A}_I=\cap_{j\in I}\mathfrak{A}_j$, and
$V_I=q^{-1}(\mathfrak{A}_I)$. Then $S_I=G_I.\exp(\mathfrak{A}_I)$
is a slice for the conjugation action of $G$ and therefore
$$G\times_{G_{I}}S_I=V_I.$$ We denote the projection to the base
by $$\pi_I: V_I\rightarrow G/G_I.$$\begin{lem}$\eta_G$ is exact
over each of the open subsets $V_j$.\end{lem}\begin{proof} $S_j':=
G_j\cdot(\mathfrak{A}_j-\mu_j)$ is a star-shaped open neighborhood
of 0 in $\mathfrak{g}_j$ and is $G_j$-equivariantly diffeomorphic
with $S_j$. We can extend this retraction from $S_j$ onto
$\exp(\mu_j)$ to a $G$-equivariant retraction from $V_j$ onto
$\mathcal{C}_j=q^{-1}(\mu_j)$. But since
$d_G\omega_{\mathcal{C}_j}+\iota^*_{\mathcal{C}_j}\eta_G=0$, then
$\eta_G$ is exact over $V_j$.
\end{proof}
Let $\iota_j :\mathcal{C}_j\rightarrow V_j$ and $\pi_j
:V_j\rightarrow G/G_j=\mathcal{C}_j$ denote the inclusion and the
projection respectively. The retraction from $V_j$ onto
$\mathcal{C}_j$ defines a $G$-equivariant homotopy operator
$$h_j :\Omega^p(V_j)\rightarrow \Omega^{p-1}(V_j).$$ Thus,
$$d_Gh_j+h_jd_G=Id-\pi_j^*\iota_j^*.$$Define the equivariant 2-form $\varpi_j$ on $V_j$ by
$(\varpi_j)_G=h_j\eta_G-\pi_j^*\omega_{\mathcal{C}_j}$.Write
$(\varpi_j)_G=\varpi_j-\theta_j$ where $\varpi_j\in \Omega^2(V_j)$
and $\theta_j\in \Omega^0(V_j,\mathfrak{g})$. For any conjugacy
class $\mathcal{C}\subset V_j$,
$\iota^*(\varpi_j)_G+\omega_{\mathcal{C}}$ is an equivariantly
closed 2-form with $\theta_j$ as its moment map. Therefore
$\iota^*(\varpi_j)_G+\omega_{\mathcal{C}}=\theta_j^*(\omega_{\mathcal{O}})_G$,
where $(\omega_{\mathcal{O}})_G$ is the symplectic form on the
(co)-adjoint orbit
$\mathcal{O}=\theta_j(\mathcal{C})$.\begin{prop}Over
$V_{ij}=V_i\cap V_j$, $\theta_i-\theta_j$ takes values in the
adjoint orbit $\mathcal{O}_{ij}$ through $\mu_i-\mu_j$.
Furthermore,
$$(\varpi_i)_G-(\varpi_j)_G=\theta_{ij}^*(\omega_{\mathcal{O}_{ij}})_G$$where
$\theta_{ij} :=\theta_i-\theta_j : V_{ij}\rightarrow
\mathcal{O}_{ij}$, and $(\omega_{\mathcal{O}_{ij}})_G$ is the
equivariant symplectic form on the orbit.\end{prop}\begin{proof}
Let $\nu : \mathfrak{A}_j\rightarrow \mathfrak{t}$ be the
inclusion map. Then $$\widetilde{h}_j\circ
(\exp\mid_{\mathfrak{A}_j})\frac{1}{2}(\theta^L+\theta^R)=\widetilde{h}_j\circ
d\nu=\nu-\mu_j$$ where $\widetilde{h}_j$ is the homotopy operator
for the linear retraction of $\mathfrak{t}$ onto $\mu_j$. This
proves that $(\exp\mid_{\mathfrak{A}_j})^*\theta_j=\nu-\mu_j$.
Therefore, for $\xi\in \mathfrak{A}_{ij}$ we have,
$$\theta_{ij}(\exp\xi)=(\xi-\mu_i)-(\xi-\mu_j)=\mu_j-\mu_i.$$Therefore
$\theta_{ij}$ takes values in the adjoint orbit through
$\mu_j-\mu_i$ by equivariance. The difference $\varpi_i-\varpi_j$
vanishes on $T$ and is therefore determined by its contractions
with generating vector fields. But $\theta_{ij}$ is a moment map
for $\varpi_i-\varpi_j$, hence $\varpi_i-\varpi_j$ equals to the
pull-back of the symplectic form on
$\mathcal{O}_{ij}$.\end{proof}\subsection{Construction of the
Basic Gerbe}\label{basic gerbe}Let $G$ be a compact, simple,
simply connected Lie group, and $B$ is an invariant inner product
 at integral level $k>0$. Use $B$ to identify
 $\mathfrak{g}\cong \mathfrak{g}^*$ and  $\mathfrak{t}\cong
 \mathfrak{t}^*$.
 We assume that under this identification, all vertices of the alcove are contained in
the weight lattice $\Lambda^*\subset\mathfrak{t}$. This is
automatic if $G$ is the special unitary group $A_d=SU(d+1)$ or the
compact symplectic group $C_d=Sp(2d)$. In general, the following
table lists the smallest integer $k$ with this
property~\cite{MR39:1590}: \\[5.pt]
\begin{center}\begin{tabular}{c|c|c|c|c|c|c|c|c|c}
$G$& $A_d$& $B_d$& $C_d$& $D_d$& $E_6$& $E_7$& $E_8$& $F_4$& $G_2$\\
\hline
$k$& 1& 2& 1& 2& 3& 12& 60& 6& 2\\
\end{tabular}
\\[40.pt]\end{center}
For constructing the basic gerbe over $G$, we pick the standard
open cover of $G$, $\mathcal{V}=\{V_i, i=0,\cdots,d\}$. For any
$\mu\in \Lambda^*$, with stabilizer $G_{\mu}$, define a line
bundle
$$L_{\mu}=G\times_{G_{\mu}} \mathbb{C}_{\mu}$$ with the unique
left invariant connection $\nabla$, where $\mathbb{C}_{\mu}$ is
the 1-dimensional $G_{\mu}$-representation with infinitesimal
character $\mu$. $L_{\mu}$ is a $G$-equivariant pre-quantum line
bundle for the orbit $\mathcal{O}=G\cdot\mu$.
Therefore$$\frac{i}{2\pi}curve_G(\nabla)=(\omega_{\mathcal{C}})_G
:=\omega_{\mathcal{O}}-\Phi_{\mathcal{O}}$$where
$\omega_{\mathcal{O}}$ is a symplectic form for the inclusion map
$\Phi_{\mathcal{O}} :\mathcal{O}\rightarrow \mathfrak{g}^*$.
Define line bundles $$L_{ij}
:=\theta_{ij}^*(L_{\mu_j-\mu_i})$$equipped with the pull-back
connection. In three fold intersection $V_{ijk}$, the tensor
product $(\delta L)_{ijk}=L_{jk}L_{ik}^{-1}L_{ij}$ is the
pull-back of the line bundle over $G/G_{ijk}$ which is defined by
the zero weight $$(\mu_k-\mu_j)-(\mu_k-\mu_i)+(\mu_j-\mu_i)=0$$of
$G_{ijk}$. Therefore it is canonically trivial with trivial
connection. The trivial sections $t_{ijk}=1$, satisfy $\delta t=1$
and $(\delta\nabla)t=0$. Define $(F_j)_G=(\varpi_j)_G$.
Since$$\delta(F)_G=\theta_{ij}^*(\omega_{\mathcal{O}_{ij}})_G=\frac{1}{(2\pi
\sqrt{-1})}curve_G(\nabla^{ij}),$$ then $\g=(\mathcal{V},L,t)$ is
a gerbe with connection $(\nabla,\varpi)$. The construction of the
basic gerbe is discussed in more general cases
in~\cite{EM,MR1968268}.
 \subsection{The
Basic Gerbe Over SU(n)} In this Section, I show that the our
construction of the basic gerbe over $SU(n)$ is equivalent to the
construction of the basic gerbe in
Gawedzki-Reis~\cite{MR2003m:81222}.\newline \indent The
\emph{special unitary group} is the classical group:$$SU(n)=\{A\in
U(n)\mid \det A=1\},$$which is a compact connected Lie group of
dimension equal to $n^2-1$ with Lie algebra equal to the
space:$$\mathfrak{su}(n)=\{A\in
L_{\mathbb{C}}(\mathbb{C}^n,\mathbb{C}^n)\mid A^{\ast}+A=0
\,\mbox{and}\, tr A=0\}.$$Any matrix $A\in SU(n)$ is conjugate to
a diagonal matrix with entries
$$\operatorname{diag}\big(\exp((2\pi
\sqrt{-1})\lambda_1(A)),\cdots, \exp((2\pi
\sqrt{-1})\lambda_n(A))\big)$$ where
$\lambda_1(A),\cdots,\lambda_n(A)\in \mathbb{R}$ are normalized by
the identity $\Sigma_{i=1}^{n}\lambda_i(A)=0$ and
\begin{equation}\lambda_1(A)\geq \lambda_2(A)\geq \cdots\geq
\lambda_n(A)\geq\lambda_1(A)-1.\end{equation}Consider the
following maximal torus of $SU(n)$,$$T=\{A\in SU(n)\mid \mbox{$A$
is diagonal}\}.$$Let $\mathfrak{t}$ be the Lie algebra of $T$.
Thus,
 $\mathfrak{t}\cong\{\lambda\in\mathbb{R}^n\mid\Sigma_{i=1}^n\lambda_i=0\}.$
The roots $\alpha\in \mathcal{R}\subset \mathfrak{t}^*$ are the
linear maps:$$\alpha_{ij}:\mathfrak{t}\rightarrow
\mathbb{R},\,(\lambda_1,\cdots,\lambda_n)\mapsto\lambda_i-\lambda_j,\,i\neq
j,$$and the set of simple roots is
$$\mathcal{S}=\{\alpha_{1,2},\alpha_{2,3},\cdots,\alpha_{n-1,n}\}.$$The
lowest root is $\alpha_{n,1}$(~\cite{MR2003c:22001}, Appendix C).
Choose the following Weyl
chamber$$\mathfrak{t}_+=\{(\lambda_1,\cdots,\lambda_n)\in
\mathfrak{t}\mid
\lambda_1\geq\lambda_2\geq\cdots\geq\lambda_n\}.$$ In that case
the fundamental alcove is
$$\mathfrak{A}=\{(\lambda_1,\cdots,\lambda_n)\in \mathfrak{t}\mid
\lambda_1\geq\lambda_2\geq\cdots\geq\lambda_n\geq\lambda_1-1\}.$$The
basic inner product on $\mathfrak{t}$ is induced from the standard
basic inner product on $\mathbb{R}^n$. We can use this inner
product to identify $\mathfrak{t}\cong\mathfrak{t}^*$. Under this
identification $\alpha_{i,j}=e_i-e_j$, where $\{e_i\}_{i=1}^n$ is
the standard basis for $\mathbb{R}^n$. The fundamental weights are
given by $$\mu_i=\{\lambda\in \mathfrak{A}\mid
\lambda_1=\lambda_2=\cdots=\lambda_i>\lambda_{i+1}=\cdots
=\lambda_n=\lambda_1-1\}.$$
$$SU(n)^{reg}=\{A\in SU(n)\mid \mbox{all eigenvalues of $A$ have multiplicity one}\}\cong G
\times_T \,\operatorname{int}\mathfrak{A}\cong G/T.$$ For $i\in
\{1,\cdots,n\}$, define
$$\mathfrak{A}_i:= \{\lambda\in \mathfrak{A}\mid
\lambda_1\geq\cdots \geq\lambda_i>\lambda_{i+1}\geq\cdots
\geq\lambda_n\geq\lambda_1-1\}.$$Thus, the standard open cover for
$SU(n)$ is $\mathcal{V}=\{V_i\}_{i=1}^n$, where
$V_i=q^{-1}(\mathfrak{A}_i)$. Each $SU(n)_{ij}$ is isomorphic to
$U(n-1)$ with the center isomorphic to $U(1)$. Over the set of
regular elements all the inequalities are strict and we have n
equivariant line bundles $E_1,\cdots,E_n$ defined by the
eigenlines for the eigenvalues $\exp((2\pi
\sqrt{-1})\lambda_i(A))$. For $i<j$, the tensor product
$E_{i+1}\otimes\cdots \otimes E_j\rightarrow SU(n)^{reg}$ extends
to a line bundle $E_{ij}\rightarrow V_{ij}$. For $i<j<k$, we have
a canonical isomorphism $E_{ij}\otimes E_{jk}\cong E_{ik}$ over
$V_{ijk}$. These line bundles together with corresponding
isomorphisms define a gerbe over $SU(n)$, in Gawedzki-Reis sense,
which represents the generator of $H^3(SU(n),\mathbb{Z})$. Each
$E_i=G\times_T \mathbb{C}_{\nu_i}$ for some $\nu\in \Lambda^*$. In
fact, by using the standard action of $T\subset SU(n)$ on
$\mathbb{C}^n$, one can see that
$$\nu_i=e_i-\frac{1}{n}(1,\cdots,1).$$ Since $\mu_i=\Sigma_{k=1}^i e_k-\frac{i}{n}
\Sigma_{k=1}^n e_k,$ therefore, $\mu_i=\Sigma_{k=1}^i \nu_k.$
Recall that in Section \ref{basic gerbe}, to construct the basic
gerbe, we defined $L_{ij} :=\theta_{ij}^*(L_{\mu_j-\mu_i})$ on
$V_{ij}$. Thus, for $i<j$ we have \begin{eqnarray*}L_{ij}
&=&\theta_{ij}^*(L_{\mu_j-\mu_i})\\&=&\theta_{ij}^*(L_{\Sigma_{k=1}^j
\nu_k-\Sigma_{k=1}^i \nu_k})\\&=&\theta_{ij}^*(L_{\Sigma_{k=i+1}^j
\nu_k})=E_{ij}\end{eqnarray*} \section{The Relative Gerbe for
$\operatorname{Hol}:\label{rhol} \mathcal{A}_G(S^1)\rightarrow
G$}Denote the affine space of connection on the trivial bundle
$S^1\times G$ by $\mathcal{A}_G(S^1)$. Thus,
$\mathcal{A}_G(S^1)=\Omega^1(S^1,\mathfrak{g})$. The loop group
$LG=Map(S^1,G)$ acts on $\mathcal{A}_G(S^1)$ by gauge
transformations:\begin{equation}\label{hol}g\cdot
A=Ad_g(A)-g^*\theta^R.\end{equation}Taking the holonomy of a
connection defines a smooth map $$\operatorname{Hol}:
\mathcal{A}_G(S^1)\rightarrow G$$with equivariance property
$\operatorname{Hol}(g\cdot A)=Ad_{g(0)}\operatorname{Hol}(A)$. If
$\mathfrak{g}$ carries an invariant inner product $B$, we write
$L\mathfrak{g}^*$ instead of $\mathcal{A}_G(S^1)$ using the
natural pairing between $\Omega^1(S^1,\mathfrak{g})$ and
$L\mathfrak{g}=\Omega^0(S^1,\mathfrak{g})$. We will refer to this
action as the coadjoint action. However, notice that the action
\ref{hol} is not the point wise action. Recall that we have
constructed in Section \ref{basic gerbe},
\\ \indent a) An open cover $\mathcal{V}=\{V_0,\cdots,V_d\}$ of
$G$ such that $V_j/G=\mathfrak{A}_j$, where $d=rank G$.\\
\indent b) For each $V_j$, a unique $G$-equivariant deformation
retraction on to a conjugacy class
$\mathcal{C}_j=G\cdot\exp(\mu_j)$, where $\mu_j$ is the vertex of
$\mathfrak{A}_j$. This deformation retraction descends to the
linear retraction of $\mathfrak{A}_j$ to $\mu_j$.\\ \indent c)
2-forms $\varpi_j\in \Omega^2(V_j)$, with
$d\varpi_j=\eta\mid_{V_j}$, such that the pull-back onto
$\mathcal{C}_j$ is the invariant 2-form for the conjugacy class
$\mathcal{C}_j$.\begin{lem}There exists a unique $LG$-equivariant
retraction from $\widetilde{V}_j:=\operatorname{Hol}^{-1}(V_j)$
onto the coadjoint orbit $\mathcal{O}_j=LG\cdot\mu_j$, descending
to the linear retraction of $\mathfrak{A}_j$ onto
$\mu_j$.\end{lem}\begin{proof}The holonomy map sets up a
one-to-one correspondence between the sets of $G$-conjugacy
classes and coadjoint $LG$-orbits, hence both are parameterized by
points in the alcove. The evaluation map $LG\rightarrow
G,\,g\mapsto g(1)$ restricts to an isomorphism $(LG)_j\cong G_j$
~\cite{MR2002g:53151}.
Hence,$$LG\times_{(LG)_j}\widetilde{S_j}=\widetilde{V_j}$$where
$\widetilde{S_j}=\operatorname{Hol}^{-1}(S_j)$ and
$S_j=G_j\cdot\exp(\mathfrak{A}_j)$. Therefore the unique
equivariant retraction from $V_j$ onto $\mathcal{C}_j$, which
descends to the linear retraction of $\mathfrak{A}_j$ onto the
vertex $\mu_j$, lifts to the desired retraction.
\end{proof}Consider the following commutative diagram:\\*\[\begin{CD}
   \widetilde{V}_j  @>\widetilde{\pi}_j>>     \mathcal{O}_j\\
   @V\operatorname{Hol}VV   @VV\operatorname{Hol}V\\
   V_j@>\pi_j>>\mathcal{C}_j
   \end{CD}\]\\*where $\widetilde{\pi}_j:\widetilde{V}_j
   \rightarrow \mathcal{O}_j$ is the projection that we get from
   retraction. Let $\sigma_j\in \Omega^2(\widetilde{V}_j)$ denote
   the pull-backs under $\widetilde{\pi}_j$ of the symplectic forms
   on $\mathcal{O}_j$.\begin{lem}On overlaps
   $\widetilde{V}_j\cap\widetilde{V}_{j'}$,
   $\sigma_j-\sigma_{j'}=\operatorname{Hol}^*(\varpi_j-\varpi_{j'})$.\end{lem}
   \begin{proof}Both sides are closed $LG$-invariant forms, for
   which the pull-back to $\mathfrak{t}\subset L\mathfrak{g}^*$
   vanishes. Hence, it suffices to check that at any point $\mu\in
   \mathcal{O}_j\cap \mathcal{O}_{j'}\subset \mathfrak{t}\subset
   L\mathfrak{g}^*$, the contraction with $\zeta_{L\mathfrak{g}^*}$
   are equal for $\zeta\in L\mathfrak{g}$. We have \begin{eqnarray*}
   \iota(\zeta)\sigma_j&=&\widetilde{\pi}_j^*dB(\Phi_j,\zeta)\\
   &=& dB(\widetilde{\pi}_j^*\Phi_j,\zeta)\end{eqnarray*}
   where $\Phi_j:\mathcal{O}_j\hookrightarrow L\mathfrak{g}^*$ is inclusion.
   But\begin{eqnarray*}(\widetilde{\pi}_j^*\Phi_j-\widetilde{\pi}_{j'}^*\Phi_{j'})
   \mid_{g\cdot\mu}&=&g\cdot\mu_j-g\cdot\mu_{j'}\\&=&(Ad_g(\mu_j)-g^*\theta^R)
   -(Ad_g(\mu_{j'})-g^*\theta^R)\\&=&
   Ad_g(\mu_j-\mu_{j'}).\end{eqnarray*}This, however, is another
   moment map for
   $\operatorname{Hol}^*(\varpi_j-\varpi_{j'})$.\end{proof}
   We conclude that the locally defined forms
   $$\varpi\mid_{\widetilde{V}_j}=\operatorname{Hol}^*(\varpi_j)-\sigma_j$$
   patch together to define a global 2-form $\varpi\in
   \Omega^2(L\mathfrak{g}^*)$. Form the properties of $\varpi_j$ and
   $\sigma_j$ we read off:\\ \indent(i)
   $d\varpi=\operatorname{Hol}^*\eta$,\\ \indent(ii)$\iota(\zeta_{L\mathfrak{g}^*})\varpi=
   \frac{1}{2}B(\operatorname{Hol}^*(\theta^L+\theta^R),\zeta(0))-dB(
   \mu,\zeta)$. \\ Here $\mu:L\mathfrak{g}^*\rightarrow
   L\mathfrak{g}^*$ is the identity map. Such a 2-form was
   constructed in~\cite{MR99k:58062} using a different method. \\ \indent
   Consider the case $G=SU(n)$ or $G=Sp(n)$, which are the two
   cases where the vertices of the alcove lie in the weight
   lattice $\Lambda^*$, where we identify $\mathfrak{g}^*\cong
   \mathfrak{g}$ using the basic inner product.
   Let$$U(1)\rightarrow \widehat{LG}\rightarrow LG$$ denote the
  $k$-th power of the basic central extension of the loop group~\cite{MR88i:22049}.
  That is, on the Lie algebra level the central
  extension$$\mathbb{R}\rightarrow\widehat{L\mathfrak{g}}\rightarrow
  L\mathfrak{g}$$ is defined by the cocycle,$$(\xi_1,\xi_2)\mapsto
  \int_{S^1}B(\xi_1,d\xi_2)\qquad \xi\in
  L\mathfrak{g}=\Omega^0(S^1,\mathfrak{g}))$$where $B$ is the
  inner product at level $k$. The coadjoint action of
  $\widehat{LG}$ on $\widehat{L\mathfrak{g}^*}=L\mathfrak{g}^*\times \mathbb{R}$
  preserves the level sets $L\mathfrak{g}^*\times\{t\}$,
  and the action for $t=1$ is exactly the gauge action of $LG$
  considered above.
   Since $\mu_j\in \Lambda^*$, the orbits $LG\cdot
   \mu_j=\mathcal{O}_j$ carry $\widehat{LG}$-equivariant
   pre-quantum line bundles
   $L_{\mathcal{O}_j}\rightarrow \mathcal{O}_j$, given explicitly
   as
   $$L_{\mathcal{O}_j}=\widehat{LG}\times_{(\widehat{LG})_j}\mathbb{C}_{(\mu_j,1)}.$$
   Here $(\widehat{LG})_j$ is the restriction of $\widehat{LG}$ to
   the stabilizer $(LG)_j$ of $\mu_j \in \mathfrak{t}\subset
   L\mathfrak{g}^*$, and $\mathbb{C}_{(\mu_j,1)}$ denotes the
   1-dimensional representation of $(\widehat{LG})_j$ with weight
   $(\mu_j,1)\in \Lambda^*\times \mathbb{Z}$. Let $E_j\rightarrow
   \widetilde{V}_j$ be the pull-back
   $\widetilde{\pi}_j^*L_{\mathcal{O}_j}$. On overlaps,
   $\widetilde{V}_j\cap \widetilde{V}_{j'}$, $E_j\otimes E_{j'}^{-1}$ is an
   associated bundle for the weight
   $(\mu_j,1)-(\mu_{j'},1)=(\mu_j-\mu_{j'},0)$. Therefore, $E_j\otimes
   E_{j'}^{-1}$ is an $LG$-equivariant bundle. It is clear by
   construction that $E_j\otimes E_{j'}^{-1}$ is the pull-back of the
   pre-quantum line bundle over $\mathcal{O}_{jj'}\subset
   \mathfrak{g}^*$. Taking all this information together, we have
   constructed an explicit quasi-line bundle for the pull-back of
   the $k$-th power of the basic gerbe under holonomy map
   $\operatorname{Hol}:L\mathfrak{g}^*\rightarrow G$ with error
   2-form equal to $\varpi\in \Omega^2(L\mathfrak{g}^*)$.
\section{Review of Group-Valued Moment Maps}
Suppose $(M,\omega)$ is a symplectic manifold together with a
symplectic action of a Lie group $G$. This action called
Hamiltonian if there exists a smooth equivariant map
\[\Phi: M\rightarrow \mathfrak{g}^*\]such that \[\iota
(\xi_M)\omega+d\langle\Phi,\xi\rangle=0\]for all $\xi\in
\mathfrak{g}$, where $\xi_M$ is the vector field on $M$ generated
by $\xi\in \mathfrak{g}$,
i.e.,\[\xi_M(m)=\frac{d}{dt}\mid_{t=0}\exp(t\xi)\cdot m.\] The map
$\Phi$ and the triple $(M,\omega,\Phi)$ are known as moment map
and Hamiltonian $G$-manifold respectively~\cite{MR91d:58073}. Let
$G$ be a compact Lie group. Fix an invariant inner product $B$ on
$\mathfrak{g}$, which we use to identify
$\mathfrak{g}\cong\mathfrak{g}^*$. Since the exponential map
$exp:\mathfrak{g}\rightarrow G$ is a diffeomorphism in a
neighborhood of the origin, the composition map \[\Psi:=exp\circ
\Phi:M\rightarrow G\] inherits the properties of the moment map
$\Phi$ and vice versa.
\begin{defn} A quasi-Hamiltonian $G$-space with group-valued moment map
is a triple $(M,\omega,\Psi)$ consisting of a $G$-manifold M, an
invariant 2-form $\omega\in\Omega^2(M)$, and an equivariant smooth
map $\Psi:M\rightarrow G$ such that
\begin{enumerate}\item $d\omega =\Psi^{*}\eta$ where $\eta \in
\Omega^{3}(G)$ is the 3-form defined by $B$. This condition is
called the relative cocycle condition.\item
$\iota(\xi_{M})\omega=\frac{1}{2}B(\Psi^{*}(\theta^{L}+\theta^{R}),\xi)$.
This condition is called the moment map condition.\item The $\ker(
\omega_{m})\in T_m(M)$ for $m\in M$ consists of all $\xi_{M}(m)$
such that
\[(Ad_{\Psi(m)}+1)\xi=0.\]This is called the minimal degeneracy
condition.\end{enumerate}
\end{defn}\subsection{Examples}
\begin{exmp}Consider a Hamiltonian $G$-manifold $(M,\omega,\Phi)$ such that
the image of $\Phi$ is a subset of the set of regular values for
the exponential map. Then $(M,\Upsilon,\Psi)$ is a
quasi-Hamiltonian $G$-space with group-valued moment map, where
\[\Psi=\exp\circ\Phi\]and
\[\Upsilon:=\omega+\Phi^*\varpi\]where $\varpi\in \Omega^2(\mathfrak{g})$
is the primitive for $\exp^*\eta$ given by the de Rham homotopy
operator for the vector space $\mathfrak{g}$. The converse is also
true, provided that $\Psi(M)$ lies in a neighborhood of the origin
on which the exponential map is a
diffeomorphism.\end{exmp}\begin{exmp}\label{Lisa} Let
$\mathcal{C}\subset G$ be a conjugacy class of $G$. The triple
$(\mathcal{C},\omega,\Phi)$ is a quasi-Hamiltonian $G$-space with
group-valued moment map where $\Phi:\mathcal{C}\hookrightarrow G$
is inclusion and
$\omega_g(\xi_{\mathcal{C}}(g),\zeta_{\mathcal{C}}(g))=
\frac{1}{2}B((Ad_g-Ad_{g^{-1}})\xi,\zeta)$
~\cite{MR98e:58034}.\end{exmp}\begin{exmp}Given an involutive  Lie
group automorphism $\rho\in Aut(G)$, i.e., $\rho^2=1$, one defines
\emph{twisted conjugacy classes} to be the orbits of the action
$h\cdot g=\rho(h)gh^{-1}.$ $G$ is a symmetric space
$$G=G\times G/(G\times G)^\rho$$ where $\rho(g_1,g_2)=(g_2,g_1)$.
The map $G\times G\rightarrow \mathbb{Z}_2\ltimes G\times G,
(g_1,g_2)\mapsto(\rho^{-1},g_1,g_2)$ takes the twisted conjugacy
classes of $G\times G$ to conjugacy classes of the disconnected
group $\mathbb{Z}_2\ltimes G\times G$. Thus by using example
\ref{Lisa} the group $G$ itself becomes a group-valued Hamiltonian
$\mathbb{Z}_2\ltimes G\times G$, with 2-form $\omega=0$, moment
map $g\mapsto (\rho,g,g^{-1})$ and action $(g_1,g_2)\cdot g=g_2 g
g_1^{-1}, \rho\cdot g=g^{-1}$.\end{exmp}\begin{exmp}Let $D(G)$ be
a product of two copies of $G$. On $D(G)$, we can define a
$G\times G$ action
by\[(g_1,g_2).(a,b)=(g_1ag_2^{-1},g_2ag_1^{-1}).\]Define a
map\[\Psi:D(G)\rightarrow G\times G,\quad
\Psi(a,b)=(ab,a^{-1}b^{-1})\]and let the 2-form $\omega$ be
defined by
\[\omega=\frac{1}{2}(B(\operatorname{Pr}_1^{*}\theta^L,
\operatorname{Pr}_2^{*}\theta^R)+B(\operatorname{Pr}_1^{*}\theta^R,\operatorname{Pr}_2^{*}\theta^L))\]where
$\operatorname{Pr}_1$ and $\operatorname{Pr}_2$ are projections to
the first and second factor. Then the triple $(D(G),\omega,\Psi)$
is a Hamiltonian $G\times G$-manifold with group-valued moment
map.\end{exmp}
\begin{exmp}\label{examplelisa}
Let $G=SU(2)$ and $M=S^4$ the unit sphere in
$\mathbb{R}^5\cong\mathbb{C}^2\times\mathbb{R}$, with SU(2)-action
induced from the action on $\mathbb{C}^2$. $M$ carries the
structure of a group-valued Hamiltonian $SU(2)$-manifold, with the
moment map $\Psi:M\rightarrow SU(2)\cong S^3$ the suspension of
the Hopf fiberation $S^3\rightarrow S^2.$ For details,
see~\cite{MR2003d:53151}. This example is generalized by
Hurtubise-Jeffrey-Sjamaar in~\cite{2} to $G=SU(n)$ acting on
$M=S^{2n}$ (viewed as unit sphere in $\mathbb{C}^n\times
\mathbb{R})$.
\end{exmp}
The equivariant de
Rham complex is defined as
\[\Omega^{k}_G(M)= \bigoplus_{2i+j=k}(\,\Omega^{j}(M)\otimes S^i(\mathfrak{g}^*)
)^G\]where $S(\mathfrak{g}^*)$ is the symmetric algebra over the
dual of the Lie algebra of $G$. Elements in this complex can be
viewed as equivariant polynomial maps from $\mathfrak{g}$ into the
space of differential forms. $\Omega_G(M)$ carries an equivariant
differential $d_G$ of degree 1,
\[(d_G\alpha)(\xi):=d\alpha(\xi)+\iota(\xi_M)\alpha(\xi).\]Since
$(d+\iota(\xi_M))^2=L(\xi_M)$ and we are restricting on the
equivariant maps, $d_G^2=0$. The equivariant cohomology is the
cohomology of this co-chain complex~\cite{MR2001i:53140}. The
canonical 3-form $\eta$ has a closed equivariant extension
$\eta_G\in \Omega_G^3(G)$ given by
$$\eta_G(\xi):=\eta+\frac{1}{2}B(\theta^L+\theta^R,\xi).$$
We can combine the first two conditions of the definition of a
group-valued moment map and get the condition
$$d_G\omega=\Psi^*\eta_G.$$
\subsection{Products}Suppose $(M,\omega, (\Psi_1,\Psi_2))$ is a
group-valued Hamiltonian $G\times G$-manifold. Then
$\widetilde{M}=M$ with diagonal action, moment map
$\widetilde{\Psi}=\Psi_1\Psi_2$ and 2-form
$$\widetilde{\omega}=\omega-\frac{1}{2}
B(\Psi_1^{*}\theta^L,\Psi^*_{2}\theta^R)$$ is a group-valued
quasi-Hamiltonian $G$-space. If $\widetilde{M}=M_1\times M_2$ is a
direct product of two group-valued quasi-Hamiltonian $G$-spaces,
we call $\widetilde{M}$ the fusion product of $M_1$ and $M_2$.
This product is denoted by $M_1\circledast M_2$. If we apply
fusion to the double $D(G)$, we obtain a group-valued
quasi-Hamiltonian $G$-space with $G$-action
$$g\cdot(a,b)=(Ad_ga,Ad_gb),$$moment map
$$\Psi(a,b)=aba^{-1}b^{-1}\equiv[a,b],$$ and 2-form
$$\omega=\frac{1}{2}(B(\operatorname{Pr}_1^{*}\theta^L,\operatorname{Pr}_2^{*}\theta^R)
+B(\operatorname{Pr}_1^{*}\theta^R,\operatorname{Pr}_2^{*}\theta^L)-
B((ab)^{*}\theta^L,(a^{-1}b^{-1})^*\theta^R)).$$Fusion of $h$
copies of $D(G)$ and conjugacy classes
$\mathcal{C}_1,\cdots,\mathcal{C}_r$ gives a new quasi-Hamiltonian
space with the moment
map\[\Psi(a_1,b_1,\cdots,a_h,b_h,d_1,\cdots,d_r)=\prod^h_{j=1}[a_j,b_j]\prod^r_{k=1}d_k.\]
\subsection{Reduction}The symplectic reduction works as usual:
\\ \indent If $(M,\omega,\Psi)$ be a Hamiltonian $G$-space with group-valued moment map
and the identity element $e\in G$ be a regular value of $\Psi$,
then $G$ acts locally freely on $\Psi^{-1}(e)$ and therefore
$\Psi^{-1}(e)/G$ is smooth. Furthermore, the pull- back of
$\omega$ to identity level set descends to a symplectic form on
$M//G:=\Psi^{-1}(e)/G$. For instance, the moduli space of flat
G-bundles on a closed oriented surface of genus h with $r$
boundary components, can be written
\[\mathcal{M}(\Sigma;\mathcal{C}_1,\cdots,\mathcal{C}_r)=G^{2h}\circledast\mathcal{C}_1
\circledast\cdots\circledast\mathcal{C}_r//G=
\Psi^{-1}(e)/G\]where the $j$-th boundary component is the bundle
corresponding to the conjugacy class $\mathcal{C}_j$. More details
can be found in ~\cite{MR99k:58062},~\cite{MR2003d:53151}.

\section{Pre-quantization of $G$-Valued Moment Maps}
We know that a symplectic manifold $(M,\omega)$ is pre-quantizable
(admits a line bundle L over M with curvature 2-form $\omega$) if
the 2-form $\omega$ is integral. In this Section, we will first
introduce a notion of a pre-quantization of a space with G-valued
moment map and then  give a similar criterion for being
pre-quantizable.\begin{defn} Let G be a compact connected Lie
group with canonical 3-form $\eta$. Fix a gerbe $\mathcal{G}$ on
$G$ with connection $(\nabla,\varpi)$ such that
$curv(\mathcal{G})=\eta$. A pre-quantization of $(M,\omega,\Psi)$
is a relative gerbe with connection $(\mathcal{L},\g)$
corresponding to the map $\Psi$ with relative curvature
$(\omega,\eta)$.
\end{defn}Since $\eta$ is closed 3-form and $\Psi^*\eta=d\omega$,
$(\omega,\eta)$ defines a relative cocycle. Recall from chapter 1
that a class $[(\omega,\eta)]\in H^3(\Psi,\mathbb{R})$ is integral
if and only if $\int_{\beta}\eta-\int_\Sigma\omega\in \mathbb{Z}$
for all relative cycles $(\beta,\Sigma)\in C_3(\Psi,\mathbb{R})$.
\begin{rem}
\label{integral} $(M,\omega,\Psi)$ is pre-quantizable if and only
if $[(\omega,\eta)]$ is integral by Theorem \ref{quantization}.
\end{rem}
\begin{thm}Suppose $M_i$, $i=1,2$ are two
quasi-Hamiltonian $G$-spaces. The fusion product $M_1\circledast
M_2$ is pre-quantizable if both $M_1$ and $M_2$ are
pre-quantizable.\end{thm}\begin{proof}Let $\operatorname{Mult}:
G\times G\rightarrow G$ be group multiplication and
$\operatorname{Pr}_i: G\times G\rightarrow G$, $i=1,2$ projections
to the first and second factors.
   Since $$\operatorname{Mult}^*\eta=\operatorname{Pr}_1^*\eta+\operatorname{Pr}_2^*\eta
   +\frac{1}{2}B(\operatorname{Pr}_1^*\theta^L,\operatorname{Pr}_1^*\theta^R),$$
   we get a quasi-line bundle with connection for the gerbe
   $\operatorname{Mult}^*\g\otimes (\operatorname{Pr}_1^*\g)^{-1}\otimes
   (\operatorname{Pr}_2^*\g)^{-1}$ such that the error 2-form is
   equal to
   $\frac{1}{2}B(\operatorname{Pr}_1^*\theta^L,\operatorname{Pr}_1^*\theta^R)$.
   Any two such quasi-line bundles differ by a flat line bundle
   with connection. Let $\Psi_i$, $i=1,2$ be moment maps
   for $M_i$, $i=1,2$ respectively and $\Psi=\Psi_1\Psi_2$
   be the moment map for their fusion product $M_1\circledast
M_2$. Thus,\begin{eqnarray*}\Psi^*\g &=&
   (\Psi_1\times \Psi_2)^*\operatorname{Mult}^*\g\\&=& (\Psi_1^*\times \Psi_2^*)
   \big((\operatorname{Pr}_1^*\g)\otimes
   (\operatorname{Pr}_2^*\g)\big)\\&=&\Psi_1^*\g\otimes
   \Psi_2^*\g.\end{eqnarray*}Therefore $M_1\circledast
M_2$ is pre-quantizable if and only if both $M_1$ and $M_2$ are
pre-quantizable.\end{proof}
\begin{prop}Suppose $G$ is simple and simply connected. Let
$k\in\mathbb{Z}$ be the level of $(M,\omega,\Psi)$. Suppose
$H^2(M,\mathbb{Z})=0$. Then there exists a pre-quantization of
$(M,\omega,\Psi)$ if and only if the image of
$$\Psi^*:H^3(G,\mathbb{Z})\rightarrow
H^3(M,\mathbb{Z})$$ is $k$-torsion.\end{prop}
\begin{proof}By assumption, $[\eta]$ represents $k$ times the generator of
$H^3(G,\mathbb{Z})$. If $H^2(M,\mathbb{Z})=0$, the long exact
sequence:\[\cdots\rightarrow H^2(G,\mathbb{Z})\rightarrow
H^2(M,\mathbb{Z})\rightarrow H^3(\Psi,\mathbb{Z})\rightarrow
H^3(G,\mathbb{Z})\overset{\Psi^*}\rightarrow
H^3(M,\mathbb{Z})\rightarrow \cdots\]shows that the map
$H^3(\Psi,\mathbb{Z})\rightarrow H^3(G,\mathbb{Z})$ is injective.
In particular, $H^3(\Psi,\mathbb{Z})$ has no torsion, and
$(M,\omega,\Psi)$ is pre-quantizable if and only if $[\eta]$ is in
the image of the map $H^3(\Psi,\mathbb{Z})\rightarrow
H^3(G,\mathbb{Z})$, i.e., in the kernel of
$H^3(G,\mathbb{Z})\rightarrow H^3(M,\mathbb{Z})$. This exactly
means that the image of this map is $k$-torsion.\end{proof}
\begin{prop}\label{condition}
If $\,H_2(M,\mathbb{Z})=0$ a pre-quantization exists. More
generally, if $H_2(M,\mathbb{Z})$ is r-torsion, a level $k$
pre-quantization exists, where $k$ is a multiple of $r$.\end{prop}
\begin{proof}If $rH_2(M,\mathbb{Z})=0$, for any
cycle $S\in C_2(M)$, there is a 3-chain $T\in C_3(M)$ with
$\partial T=r\cdot S$. If $\Psi(S)=\partial B$, $\Psi(T)-rB$ is a
cycle and
\begin{equation}\begin{split}\int_S
k\omega -\int_B k\eta &=\frac{k}{r}(\int_Td\omega-\int_{rB}\eta)\\
&=\frac{k}{r}(\int_T \Psi^*\eta-\int_{rB}\eta)\\
&=\frac{k}{r}(\int_{\Psi (T)}\eta-\int_{rB}\eta)\\
&=\frac{k}{r}(\int_{\Psi(T)-rB}\eta)\in
\mathbb{Z}.\end{split}\end{equation}By Remark \ref{integral}
$(M,\omega,\Psi)$ is
pre-quantizable.\end{proof}\begin{exmp}$M=S^{4}$ carries the
structure of a group-valued Hamiltonian $SU(2)$-manifold, with the
moment map $\Psi:M\rightarrow SU(2)\cong S^3$ the suspension of
the Hopf fiberation $S^3\rightarrow S^2$, Example
\ref{examplelisa}. By Proposition \ref{condition} this
$SU(2)$-valued moment map is pre-quantizable.\end{exmp}
\subsection{Reduction}Let $G$ be a simply connected Lie group. Fix
a pre-quantization $\lc$ for a space with G-valued moment map
$(M,\omega,\Psi)$.\[\\[7.pt]\begin{CD}
\lc  @. \g\\
@VVV   @VVV\\
M @>\Psi>>G\\
@A{\iota}AA   @AAA\\
\Psi^{-1}(e) @>\Psi>> \{e\}
\end{CD}\\[7.pt]\]
Since $\mathcal{G}|_{\Psi^{-1}(e)}$ is equal to trivial gerbe,
$\mathcal{L}|_{\Psi^{-1}(e)}$ is a line bundle with connection
with curvature $(\iota_{\Psi^{-1}(e)})^*\omega$. Since G is simply
connected and the 2-form $(\iota_{\Psi^{-1}(e)})^*\omega$ is
$G$-basic, there exists a unique lift of the $G$-action to
$\mathcal{L}|_{\Psi^{-1}(e)}$ in such a way that the generating
vector fields on $\mathcal{L}|_{\Psi^{-1}(e)}$ are horizantal.
This is a special case of Kostant's construction~\cite{MR45:3638}.
In conclusion, we get a pre-quantum line bundle over
$\Psi^{-1}(e)/G$.
\subsection{A Finite Dimensional Pre-quantum Line Bundle for
$\mathcal{M}(\Sigma)$} Let $M=G^{2h}$ where $G$ is a simply
connected Lie group and consider the map
$$\Psi:M\rightarrow G$$with the rule
\[\Psi(a_1,\cdots,a_h)= \prod_{i=1}^h[a_i,b_i].\]Let $\mathcal{G}$
be the basic gerbe with the connection on G and $curv({\g})=\eta$.
The moduli space of flat G-bundles on a closed oriented surface
$\Sigma$ of genus h is equal to
\[\mathcal{M}(\Sigma)=G^{2h}//G=\Psi^{-1}(e)/G.\]Since $G$ is simply
connected, $H^2(G,\mathbb{Z})=H^2(G^{2h},\mathbb{Z})=0$.
$H^3(G^{2h},\mathbb{Z})\simeq \mathbb{Z}$ is torsion free
therefore by Proposition 4.3.3 there exists a unique quasi-line
bundle$\mathcal{L}$ for the gerbe $\Psi^*\mathcal{G}$.
\newline \indent Pick a connection for this quasi-line bundleand call the error 2-form
$\nu$. Therefore $d(\nu-\omega)=0$. This together with the fact
that $H^2(M,\mathbb{Z})=0$ allow us to modify quasi-line bundle
with connection such that triple $(M,\omega,\Psi)$ is
pre-quantizabe. By reduction we get a pre-quantum line bundle over
$\Psi^{-1}(e)/G= \mathcal{M}(\Sigma)$.

\subsection{Pre-quantization of Conjugacy Classes of a Lie
Group}Let G be a simple, simply connected compact Lie group. Fix
an inner product $B$ at level $k$. The map
$$\exp:\mathfrak{g}\rightarrow G$$ takes
(co)adjoint orbits $\mathcal{O}_{\xi}$ to conjugacy classes
$\mathcal{C}=G\cdot{\exp(\xi)}$.\\ \indent Any conjugacy class
$\mathcal{C}\subseteq G$ is uniquely a G-valued quasi-Hamiltonian
$G$-space $(\,\mathcal{C},\omega,\Psi)$, where
$\Psi:\mathcal{C}\hookrightarrow G$ is inclusion map, as it
explained in example \ref{Lisa}. Suppose $(\beta,\Sigma)\in
\Co_n(\Psi,\mathbb{Z})$ is a cycle. We want to see under which
conditions $\mathcal{C}$ is pre-quantizable at level $k$.
Equivalently, we are looking for conditions which implies
\[k(\int_{\beta}\eta-\int_\Sigma\omega)\in \mathbb{Z}\] where
$\eta=\frac{1}{12}B(\theta^L,[\theta^L,\theta^L])$ is canonical
3-form. Consider the basic gerbe
$\mathcal{G}=(\mathcal{V},L,\theta)$ with connection on G with
curvature $\eta$. For all $\mathcal{C}\subseteq G$ there exists a
unique $\xi\in\mathfrak{A}$ such that $\exp(\xi)\in \mathcal{C}$.
Let
$$\iota_{\mathcal{C}}:\mathcal{C}\rightarrow G$$ be inclusion map
assume that $\varpi_0$ is the primitive of $\eta$ on $V_0$, i.e.,
$$\eta\mid_{V_0}=d\varpi_{0}$$where $V_{0}$ contains $\mathcal{C}$.
Recall from Section \ref{basic gerbe} that
$$\omega_{\mathcal{C}}=\theta_0^*(\omega_{\mathcal{O}_{\xi}})_G-\iota_{\mathcal{C}}^*(\varpi_0)_G$$
and pull-back of the $\theta$ to $\mathcal{C}$ is zero. Thus,
\begin{eqnarray*}k(\int_{\Sigma}\omega_{\mathcal{C}}-\int_{\beta}\iota_{\mathcal{C}}^*\eta) &=&
k(\int_{\Sigma}\theta_0^*(\omega_{\mathcal{O}_{\xi}})_G-\iota_{\mathcal{C}}^*(\varpi_0)_G-\int_{\beta}
\iota_{\mathcal{C}}^*(\eta))\\
&=&k(\int_{\Sigma}\theta_0^*(\omega_{\mathcal{O}_{\xi}})_G).\end{eqnarray*}
$(\omega_{\mathcal{C}},\iota_{\mathcal{C}}^*\eta)$ is integral if
and only if the symplectic 2-form
 $k\omega_{\mathcal{O}_{\xi}}$ is
integral. It is a well-known fact from symplectic geometry that
$k\omega_{\mathcal{O}_{\xi}}$ is integral if and only if
$B(\xi)\in \Lambda^*_k:=\Lambda^*\cap k\mathfrak{A}$, by viewing
$B$ as a linear map $\mathfrak{t}\rightarrow \mathfrak{t}^*$.
   \section{Hamiltonian Loop Group Spaces}
   Fix an invariant inner product $B$ on $\mathfrak{g}$.
   Assume that $G$ is simple and simply connected.
   Recall that a Hamiltonian loop group manifold is a triple
   $(\widetilde{M},\widetilde{\omega},\widetilde{\Psi})$ where
   $\widetilde{M}$ is an (infinite-dimensional) $LG$-manifold,
   $\widetilde{\omega}$ is an invariant symplectic form on $\widetilde{M}$,
   and $\widetilde{\Psi}:\widetilde{M}\rightarrow L\mathfrak{g}^*$
   an equivariant map satisfying the usual moment map condition,
   $$\iota(\xi_{\widetilde{M}})\widetilde{\omega}+dB(\widetilde{\Psi},\xi)=0
   \qquad\qquad\xi\in\Omega^0(S^1,\mathfrak{g}).$$
   \begin{exmp}Let $\mathcal{O}\subset
   L\mathfrak{g}^*$
   be an orbit of the loop group action. Then $\mathcal{O}$
   carries a unique structure for a Hamiltonian $LG$-manifold when
   the moment map is inclusion and the 2-form is $$
   \widetilde{\omega}_{\mu}(\xi_{\mathcal{O}}(\mu),\eta_{\mathcal{O}}(\mu))=\langle
   d_{\mu}\xi,\eta\rangle=\int_{S^1}B((d_{\mu}\xi),\eta).$$\end{exmp}
   The based loop group $\Omega G\subset LG$ consisting of loops that are trivial
   at the origin of $S^1$, acts freely on $L\mathfrak{g}^*$
   and the quotient is just the holonomy map. There
   is a one-to-one correspondence between quasi-Hamiltonian $G$-spaces
   $(M,\omega,\Psi)$ and Hamiltonian $LG$-spaces with proper
   moment maps
   $(\widetilde{M},\widetilde{\omega},\widetilde{\Psi})$, where
   \begin{eqnarray*}M &=& \widetilde{M}/\Omega G,\\
   \operatorname{Hol}\circ \widetilde{\Psi}&=& \Psi\circ
   \operatorname{Hol},\\
   \widetilde{\omega}&=&\operatorname{Hol}^*\omega-\widetilde{\Psi}\varpi.
   \end{eqnarray*}This is called \emph{Equivalence Theorem} in~\cite{MR99k:58062}.
   We thus, have a commutative diagram:
   \\*\[\begin{CD}
   \widetilde{M}  @>\widetilde{\Psi}>>     L\mathfrak{g}^*\\
   @V{\operatorname{Hol}}VV   @V{\operatorname{Hol}}VV\\
   M@>\Psi>>G
   \end{CD}\]\\\begin{thm}(Equivalence Theorem for Pre-quantization) There is a one-to-one correspondence between
   pre-quantizations of quasi-Hamiltonian $G$-spaces with group
   valued moment maps and pre-quantizations of the corresponding Hamiltonian $LG$-spaces
with proper moment maps.\end{thm}\begin{proof} Assume that we have
constructed a pre-quantization of a quasi-Hamiltonian $G$-space
with group-valued moment map $(M,\omega,\Psi)$ with the
corresponding Hamiltonian $LG$-space
$(\widetilde{M},\widetilde{\omega},\widetilde{\Psi})$. Thus, we
have a relative gerbe mapping to the basic gerbe over $G$.
Pull-back of this quasi-line bundle under
$\operatorname{Hol}:\widetilde{M}\rightarrow M$, gives a
quasi-line bundle of the gerbe
$\operatorname{Hol}^*\Psi^*\g=\widetilde{\Psi}^*\operatorname{Hol}^*\g$
over $\widetilde{M}$. But recall that, we have a quasi-line bundle
for $\operatorname{Hol}^*\g$ as it explained in Section
\ref{rhol}. Therefore the difference between these two quasi-line
bundles with connection is a line bundle with connection
$\widetilde{L}\rightarrow \widetilde{M}$ with the curvature 2-form
$\widetilde{\omega}=\operatorname{Hol}^*\omega-\widetilde{\Psi}\varpi$
by Remark \ref{qc}. Note also that if the quasi-line bundle for
$\Psi:M\rightarrow G$ is $G$-equivariant, then since the
quasi-line bundle for $L\mathfrak{g}^*\rightarrow G$ is
$\widehat{LG}$-equivariant, the line bundle $\widetilde{L}$ will
be $\widehat{LG}$ equivariant. Conversely, suppose that we are
given a $\widehat{LG}$-equivariant line bundle over
$\widetilde{M}$, where $U(1)\subset \widehat{LG}$ acts with weight
1. The difference of this $\widehat{LG}$-equivariant line bundle
and the $\widehat{LG}$-equivariant quasi-line bundle for
$\operatorname{Hol}^*\Psi^*\g$ (constructed in Section
\ref{rhol}), is a quasi-line bundle with error 2-form
$\operatorname{Hol}^*\omega$. By descending of this quasi-line
bundle to $M$, we can get the desired quasi-line bundle for
$\Psi^*\g$.
\end{proof}
The argument, given here applies in greater generality:\newline
For any $\widehat{LG}$-equivariant line bundle
$\widetilde{L}\rightarrow \widetilde{M}$, where the central
extension $U(1)\subset\widehat{LG}$ acts with weight
$k\in\mathbb{Z}$, there is a corresponding relative gerbe at level
$k$ with respect to the map $\Psi:M\rightarrow G$. Indeed, the
given quasi-line bundle for
$\widetilde{\Psi}^*\operatorname{Hol}^*\g^k$ is given by
$\widehat{LG}$ equivariant line bundles over
$\operatorname{Hol}^{-1}\Psi^{-1}(V_j)$ at level $k$. Twisting by
$\widetilde{L}$, we get new quasi-line bundle where $U(1)\subset
\widehat{LG}$ acts trivially. The quotient therefore descends to a
quasi-line bundle over $M$. For instance, Meinrenken and Woodward
construct for any Hamiltonian loop group space a so-called
``canonical line bundle'' in~\cite{MR2002g:53151}, which is
$\widehat{LG}$-equivariant at level $2c$, where $c$ is the dual
Coxeter number. Therefore this line bundle gives rise to a
distinguished element of $H^3(\Phi,\mathbb{Z})$ at level $2c$.
Notice that $M$ and $\widetilde{M}$ are not pre-quantizable
necessarily.


\bibliographystyle{amsplain}
\bibliography{xbib}
\end{document}